%This is in plain Tex.
\magnification=\magstep1
\tolerance=20000
\overfullrule=0pt
\magnification=\magstep1
\mathsurround=2pt
%\nopagenumbers
%\headline={ifodd\pageno=1\rightheadline \else \leftheadline\fi}
%\def\rightheadline{\tenrm\hfil RIGHT RUNNING HEAD\hfill\folio}
%\def\leftheadline{\tenrm\hfil LEFT RUNNING HEAD\hfill}
%\voffset=2\baselineskip
\footline={\ifnum\pageno=1{\hfil}\else{\hss\tenrm\folio\hss}\fi}

\font\tenmsbm=msbm10
\font\sevenmsbm=msbm7
\font\fivemsbm=msbm5
\newfam\msbmfam
\textfont\msbmfam=\tenmsbm
\scriptfont\msbmfam=\sevenmsbm
\scriptscriptfont\msbmfam=\fivemsbm
\def\Bbb#1{\fam\msbmfam#1}
\baselineskip=20pt
\def\tphit{{ {}_3\phi_2}}
\def\frac#1#2{{#1\over #2}}
\def\bR{{\Bbb R}}
\def\bC{{\Bbb C}}

\def\min{{\rm min}}
\def\tX{{\widetilde X}}
\def\appX{{\buildrel \approx\over X}}

\centerline{\bf Contiguous Relations, Basic Hypergeometric Functions,}
\centerline{\bf and Orthogonal Polynomials.}
\bigskip

\centerline{\bf III. Associated Continuous Dual $q$-Hahn Polynomials
\footnote*{\rm Research partially supported by NSF grant DMS 9203659 
and NSERC grant A5384.} }
\bigskip
\baselineskip=12pt
\centerline{D.P. Gupta}
\centerline{Department of Mathematics}
\centerline{University of Toronto}
\centerline{Toronto, M5S 1A1 Canada}
\bigskip
\centerline{M.E.H. Ismail}
\centerline{Department of Mathematics}
\centerline{University of Florida}
\centerline{Tampa, Florida 33620-5700}
\centerline{U.S.A.}
\bigskip
\centerline{D.R. Masson}
\centerline{Department of Mathematics}
\centerline{University of Toronto}
\centerline{Toronto, M5S 1A1 Canada}
\bigskip

\midinsert \narrower\narrower
\noindent{\bf Abstract.} Explicit solutions for 
the three-term recurrence
satisfied by associated continuous dual $q$-Hahn polynomials are
obtained.  A minimal solution is identified and an explicit expression for the
related continued fraction is derived.  The absolutely continuous
component of the spectral measure is obtained. Eleven limit cases are
discussed in some detail.  These include associated big $q$-Laguerre
, associated Wall, associated Al-Salam-Chihara, associated Al-Salam-Carlitz
 I, and associated continuous $q$-Hermite polynomials.
\endinsert
\vglue .5in

\bigskip


\vfill\eject

\noindent{\bf 1. Introduction.}

The discrete $q$-Hahn and the discrete dual $q$-Hahn polynomials were
introduced by Hahn in [10] and [11] respectively.  These are particular
cases of polynomials introduced by Askey and Wilson and called $q$-Racah
polynomials [3].  Askey and Wilson [3] have also considered 
continuous dual $q$-Hahn polynomials as a particular case of the 
Askey-Wilson polynomials.  The objective of the present study is 
to generalize the continuous dual $q$-Hahn polynomials 
to the associated continuous dual
$q$-Hahn polynomials. The $q\to1$ limit gives the case of associated 
continuous Hahn polynomials which have have been studied by Ismail et al
[13]. It may be mentioned that in two earlier
papers [7], [8], we have discussed the associated continuous Hahn
(for continuous Hahn polynomials see Askey [1])
and the associated big $q$-Jacobi polynomials respectively.
In both of these associated cases we made extensive use of contiguous
relations for hypergeometric and $q$-hypergeometric functions.
Although the use of contiguous relations in connection with continued
fractions goes back to Gauss [23], the importance of contiguous
relations in relation to the theory of orthogonal polynomials was
first stressed by Wilson [25].

In Section 2, we obtain six solutions to the 
three-term recurrence relation satisfied by associated
continuous dual $q$-Hahn polynomials.  This is done with 
the help of three-term contiguous relations satisfied
by balanced ${}_3\phi_2$'s.  It is also demonstrated how an 
existing three-term transformation formula for balanced
${}_3\phi_2$'s connects any three of these solutions.
By examining the large $n$ asymptotics of the solutions and the
associated second order difference equation we show that one of
the solutions is a minimal [6],[16] solution.  When one of 
the four parameters is equal to `$q$', another of our solutions
reduces to the continuous dual $q$-Hahn polynomial solution [3].

In Section 3, the related infinite continued fraction is obtained.
Following the procedure employed in several other cases
(see [21], [22], [23]) we then derive the explicit weight
function for the absolutely continuous component of the
spectrum.

Section 4 is devoted to obtaining a generating function and 
hence an explicit expression for the associated continuous
dual $q$-Hahn polynomials.  The method employed is the same
as by Ismail and Libis [14] for big $q$-Laguerre  polynomials.

In Section 5, we examine four limiting cases of the original recurrence relation
together with their solutions, related continued fractions and 
explicit polynomials.  The first two limits are associated big
$q$-Laguerre and associated Wall polynomials.  These are on the
${}_2\phi_1$ and ${}_1\phi_1$ levels respectively.  Two further limits
are found at the ${}_0\phi_1$ level. 

In Section 6, we consider seven additional limiting cases. These include
the associated cases for Al-Salam-Chihara, Al-Salam-Carlitz I,
and continuous $q$-Hermite polynomials.

In Section 7, we give the connection between solutions to the associated
Askey-Wilson [15],[9] and the associated continuous dual $q$-Hahn
polynomial recurrence relations.
\vfill\eject

\noindent{\bf 2.  Three-term contiguous relations and solutions.}
\medskip

The recurrence relation satisfied by associated continuous
dual $q$-Hahn polynomials can be expressed as
$$
\leqalignno{
& X_{n+1} - (z-a_n) X_n + b_n^2 X_{n-1} \ = \ 0 & (2.1)\cr
a_n & := \ a_n (z;A,B,C,D)\cr
& = \ \left( {1\over A} + {1\over B} + { 1\over C} + 
{1\over D}\right) q^n-(1+q) q^{2n-1}\cr
b_n^2 & := \ b_n^2 (z;A,B,C,D)\cr
& = \ {q\over ABCD} (1-Aq^{n-1})(1-Bq^{n-1})(1-Cq^{n-1})(1-Dq^{n-1}).
\cr}
$$
The symmetry with respect to the parameters $A,B,C,D$ 
is obvious.  With
this form of the recurrence, it is easy to take successive limits 
$A,B,C,D \to \infty$.  The limit $D \to \infty$ gives
associated big $q$-Laguerre polynomials (see Ismail and Libis
[14]) and a subsequent limit $C \to \infty$ gives associated
Wall polynomials (see Chihara [4, p. 198]).  Finally $B \to
\infty$ and then $A \to \infty$ give additional cases.

Note that in (2.1) the $a_n$ can also be expressed as
$$
\eqalign{
a_n & = \ -(\lambda_n + \mu_n ) + {1\over AB} + {q\over CD}\cr
\lambda_n & = \ (1-Aq^n)(1-Bq^n)/AB\cr
\mu_n & = \ q(1-Cq^{n-1} )(1-Dq^{n-1})/CD.\cr}
$$
This means that with a renormalization and a translation
of the coordinate $z$ we may re-express (2.1) and the 
aforementioned limits as birth and death processes
with birth and death rates $\lambda_n$ and $\mu_n$
respectively [12]. A second family with $\mu_0:=0$ should also be 
investigated. This has already been done in the more general case of
associated Askey-Wilson polynomials [15]. 

The solutions to (2.1) and its limit cases will be expressed in terms of the
basic hypergeometric functions
$$
{}_r\phi_s \left(\matrix{a_1,a_2,\cdots,a_r\cr b_1,b_2,\cdots,b_s\cr}; z\right) \ = \
\sum_{k=0}^\infty {(a_1,a_2,\cdots,a_r)_k\over (b_1,b_2,\cdots,b_s,q)_k} 
[(-1)^kq^{k(k-1)/2}]^{1+s-r}z^k,\ |z| < 1,
$$
where
$$
\eqalign{
(a)_{\infty} \ = \ \prod_{j=1}^{\infty}(1-aq^{j-1}),\ (a)_n & = \
(a)_{\infty}/(aq^n)_{\infty}, n\ \hbox{integer},\cr
(a_1,a_2,\ldots ,a_m)_n & = \ \prod\limits_{k=1}^m
(a_k)_n,\ n \ \hbox{integer or} \ \infty\ .\cr}
$$
We will use the notation
$$
\phi \ := \ \phi (a,b,c,d,e) \ = \ \tphit \left(\matrix{a,b,c\cr
d,e\cr}; {de\over abc}\right), \left| {de\over abc}\right|
\ < \ 1, \leqno(2.2)
$$
and its analytic continuation for $\left|{de\over abc}\right|
\geq 1$ given by the following transformations [5, (III.9),
(III.10), p. 241]
$$
\leqalignno{
\tphit\ \left(\matrix{a,b,c\cr d,e\cr} ; {de\over abc} \right)
\ & = \ {(e/a, de/bc)_\infty\over (e,de/abc)_\infty} \ \tphit\
\left(\matrix{a, d/ b, d/ c\cr
d, de/ bc\cr}; {e\over a} \right) & (2.2a)\cr
\tphit\ \left(\matrix{a,b,c\cr d,e\cr} ; {de\over abc} \right)
& = \ {(b, {de/ ab}, {de/ bc})_\infty\over (d,e,{de/
abc})_\infty} \ \tphit\ 
\left(\matrix{{d/ b}, {e/ b}, {de/ abc}\cr
{de/ ab}, {de/ bc}\cr};  b \right). & (2.2b)\cr}
$$
To obtain solutions to (2.1) we use $\tphit$ contiguous
relations and the usual notation
$$
\eqalign{
\phi(a\pm) & = \ \phi (aq^{\pm 1}, b,c,d,e),\cr
\phi_\pm & = \ \phi (aq^{\pm 1}, bq^{\pm 1}, cq^{\pm 1},
dq^{\pm 1}, eq^{\pm 1}).\cr}
$$
Two such contiguous relations are [18]:
$$
\phi - \phi (a+)+ {(1-b)(1-c)\over (1-d)(1-e)} {de\over
abcq} \phi_+ \ = \ 0, \leqno(2.3)
$$
$$
(1-d)(1-e) \phi + (d-a)\big(1-{e\over a}\big) \phi_+ (a- ) -
(1-a) \big(1-{de\over abcq} \big) \phi_+ \ = \ 0.
\leqno(2.4)
$$
Changing $(a,b,c,d,e) \longrightarrow (a,b /q, c/
q, d /q, e /q )$ in (2.4) and then eliminating
$\phi(a+)$ and $\phi_+$ from the  resulting equation together with 
(2.3) and (2.4), we obtain
$$
\leqalignno{
& {(1-b)(1-c)(1-{d\over a})(1-{e\over a})\over (1-d)(1-e)}
{de\over bcq} \phi_+ (a-) & (2.5)\cr
& \quad - \left[ (1-a)\big(1-{de\over abcq}\big) + a\big( 1 -
{d\over aq}\big)\big( 1-{e\over aq}\big) + {de\over abcq}
(1-b)(1-c)\right] \phi\cr
& \qquad + \big(1- {d\over q}\big)\big(1-{e\over q}\big)
\phi_- (a+) \ = \ 0.\cr}
$$
With the replacements $b = Bq^n,\ c = Cq^n, \ d /a \ =
\ Dq^n, \ e /a \ = \ Aq^n$, this becomes
$$
\leqalignno{
& \qquad {(1-Bq^n)(1-Cq^n)(1-Dq^n)(1-Aq^n)\over (1-aDq^n)(1-aAq^n)}
{DAa^2\over BCq} \phi_{n+1} & (2.6)\cr
& -\left[(1-a)(1-{DAa\over BCq})+ a(1-Dq^{n-1})(1-Aq^{n-1})
+ {DAa\over BCq} (1-Bq^n)(1-Cq^n)\right] \phi_n\cr
& \qquad + (1-Daq^{n-1})(1-Aaq^{n-1}) \phi_{n-1} \ = \
0,\cr}
$$
where
$$
\phi_n \ = \ \tphit\ \left( \matrix{a,Bq^n,Cq^n\cr
aDq^n,aAq^n\cr} ; {DAa\over BC} \right)\ .
$$
We write $z \ = \ q/aDA+a/BC$ so
that
$$
\leqalignno{
a\ & = \ {BC\over 2} \left( z \pm \ \sqrt{z^2- {4q\over
ABCD}} \right) & (2.7)\cr
&\ = \ BC \lambda_\pm , \quad {\rm say.}\cr}
$$
After renormalization, (2.6) becomes (2.1) with a solution
$$
\leqalignno{
X_n^{(1),\pm} &\ = \ X_n^{(1),\pm} (z;A,B,C,D) & (2.8)\cr
&\ = \ {(A,B,C,D)_n\over (BCD\lambda_\pm , ABC\lambda_\pm
)_n} (\lambda_\pm)^n \ \tphit\ \left(\matrix{BC\lambda_\pm,
Bq^n,Cq^n\cr BCD\lambda_\pm q^n, ABC\lambda_\pm q^n\cr} ;
AD\lambda_\pm \right).\cr}
$$
We shall show later that, with a suitable choice of square
root branch, $X_n^{(1),-}$ is a minimal solution of (2.1).
Because of the symmetry in (2.1) and the fact that the minimal solution is unique up to a constant
multiple, the parameter interchanges $A \leftrightarrow C$
or $B \leftrightarrow D$ in the above solution must yield
only an $n$ independent multiple of $X_n^{(1),-}.$  Let us verify this.
If we make the interchange $A \leftrightarrow C,\
X_n^{(1),\pm}$ changes to
$$
\xi_n^{(1),\pm} \ = \ {(A,B,C,D)_n\over (BAD\lambda_\pm ,
ABC\lambda_\pm)_n} (\lambda_\pm)^n \ \tphit 
\left(\matrix{BA\lambda_\pm , Bq^n, Aq^n\cr
BAD\lambda_\pm q^n, ABC\lambda_\pm q^n\cr} ; CD\lambda_\pm \right).
\leqno(2.9)
$$
If we apply the transformation (2.2a) to (2.9) we obtain
$$
\eqalign{
\xi_n^{(1),\pm} \ = & {(A,B,C,D)_n\over (BAD\lambda_\pm,
ABC\lambda_\pm)_n} (\lambda_\pm )^n
{(AD\lambda_\pm,BCD\lambda_\pm q^n)_\infty\over (BAD
\lambda_\pm q^n, CD\lambda_\pm)_\infty}\cr
\cr
& \quad \times \tphit\ \left( \matrix{Bq^n,Cq^n, BC\lambda_\pm\cr
ABC\lambda_\pm q^n, BCD\lambda_\pm q^n\cr} ; AD\lambda_\pm
\right)\cr}
$$
which is clearly just a constant multiple of $X_n^{(1),\pm}$.

In order to obtain a second solution of (2.1) we start from
a different three-term contiguous relation [18];
$$
\displaylines{
(2.10)\qquad \qquad
[de(a-b-c) + abc(d+e+q-a-aq)] \phi\hfill\cr
\hfill {} + (1-a)(de-abcq)
\phi(a+) + bc(d-a)(e-a)\phi (a-) \ = \ 0 .\cr} 
$$
Changing $a \to aq^n$, we can write the above as
$$
\displaylines{
[de(aq^n - b-c) + abcq^n (d+e+q-aq^n - aq^{n+1})] Y_n + \hfill\cr
\hfill {} (1-aq^n)(de-abcq^{n+1})Y_{n+1} + bc(d-aq^n)
(e-aq^n)Y_{n-1} = 0,\cr}
$$
where
$$
Y_n \ = \ \tphit\ \left(\matrix{aq^n,b,c\cr d,e\cr} ; {de\over abc} q^{-n}\right).
$$
Writing
$$
a = A,\ b = AB\lambda_+ , \ c = AB\lambda_- , \ d = {Aq\over C}, \
e = {Aq\over D}\quad {\rm and}\quad {q\over CD} ({1\over b} + {1\over c}) = z,
$$
and renormalizing we again arrive at equation (2.1) with a new
solution
$$
X_n^{(2)} \ = \ {(A,B)_n\over (AB)^n} \tphit\ \left( \matrix{Aq^n,
AB\lambda_+ , AB\lambda_-\cr Aq /C, Aq /D\cr};
{1\over B} q^{-n+1} \right). \leqno(2.11)
$$
From (2.11) we obtain additional solutions by parameter 
interchanges due to the symmetry of (2.1).  However, we find that this 
is not true for the interchanges $B \leftrightarrow C$ or $B \leftrightarrow D$.
This can be seen by applying transformation (2.2a) to $X_n^{(2)}$.
We get
$$
\leqalignno{
X_n^{(2)}&\  =\ {(A,B)_n\over (AB)^n} {(q^{-n+1} /D, Aq /B)_\infty
\over (Aq /D, q^{-n+1} /B )_\infty} \ \tphit\ 
\left( \matrix{Aq^n, AD\lambda_- , AD\lambda_+\cr
Aq /C, Aq /B\cr} ; {q^{-n+1}\over D} \right)\cr
&\  =  \ {(A,D)_n\over (AD)^n} {(Aq /B, q /D)_\infty
\over (Aq /D, q /B )_\infty} \ \tphit \ 
\left( \matrix{Aq^n, AD\lambda_- , AD\lambda_+\cr
Aq /C, Aq /B\cr} ; {q^{-n+1}\over D} \right) & (2.12) \cr}
$$
where the right side is a constant multiple of the $B \leftrightarrow
D$ interchange applied to $X_n^{(2)}$.  It is similarly seen
that $B \leftrightarrow C$ does not yield a new solution.
However the interchanges $A \leftrightarrow B, \ A
\leftrightarrow C$ and $A \leftrightarrow D$ do yield the
following new solutions:
$$
\leqalignno{
X_n^{(3)}\ & =\ {(A,B)_n\over (AB)^n} \ \tphit \ 
\left( \matrix{Bq^n, AB\lambda_+ , AB\lambda_-\cr
Bq /C, Bq /D\cr} ; {q^{-n+1}\over A} \right) , &
(2.13) \cr
X_n^{(4)}\ & =\ {(B,C)_n\over (BC)^n} \ \tphit \ 
\left( \matrix{Cq^n, BC\lambda_+ , BC\lambda_-\cr
Cq /A, Cq /D\cr} ; {q^{-n+1}\over B} \right) , &
(2.14) \cr
X_n^{(5)}\ & =\ {(B,D)_n\over (BD)^n} \ \tphit \ 
\left( \matrix{Dq^n, BD\lambda_+ , BD\lambda_-\cr
Dq /C, Dq /A\cr} ; {q^{-n+1}\over B} \right) . &
(2.15) \cr}
$$
It can be shown that the three-term transformation formula
[5, (III.33), p. 245] connects $X_n^{(1)}$ with any two of
the solutions $X_n^{(2)}, X_n^{(3)}, X_n^{(4)}$ and
$X_n^{(5)}$.  One such relation works out to be
$$
\leqalignno{
& (ABC\lambda_+ , AC\lambda_-, {q\over D}, {A\over C} )_\infty
X_n^{(1),+} - (A,A\lambda_-, AB\lambda_+ , {Cq\over
D})_\infty X_n^{(4)} & (2.16) \cr
& = \ {(C,C\lambda_-, {A\over C}, {Aq\over D}, BC\lambda_+,
CD\lambda_+ , ABD\lambda_+ )_\infty\over ({C\over A},
AD\lambda_+ , BCD\lambda_+ )_\infty} X_n^{(2)}.\cr}
$$

Another three-term contiguous relation satisfied by balanced
$\tphit$'s also yields solutions to (2.1).  The required
contiguous relation, which can be deduced from (2.3), (2.4)
and (2.10), is
$$
\leqalignno{
& {(1-a)(1-b)(1-c)\over (1-d)(1-e)} {de\over abcq} (de-abcq)
\phi_+  & (2.17) \cr 
& \quad + [abc(d+e-q) + de(1+q-a-b-c)] \phi + abcq (1-{d\over
q})(1-{e\over q})\phi_- \ = \ 0.\cr}
$$
Replacing $(a,b,c,d,e)$ by $(aq^{-n}, bq^{-n}, cq^{-n}, dq^{-n},
eq^{-n})$ and writing
$$
Z_n \ = \ \tphit\ \left(\matrix{ aq^{-n}, bq^{-n}, cq^{-n}\cr
dq^{-n}, eq^{-n}\cr} ; {de\over abc} q^n \right),
$$
we have from (2.17)
$$
\leqalignno{
\ q^{-n-1} (1-{q^{n+1}\over d})(1-{q^{n+1}\over e} ) Z_{n+1}
& + \left\{ {1\over d} + {1\over e} - {q^{n+1}\over de} +
{q^n\over abc} (q^n + q^{n+1} -a-b-c) \right\} Z_n \cr
& + {q^n\over de} {(1-{1\over a} q^n)(1 - {1\over
b}q^n) (1-{1\over c}q^n)(1-{de\over abcq} q^n)\over
(1-{1\over d} q^n )(1-{1\over e}q^n)} Z_{n-1} \ = \ 0. & (2.18) \cr}
$$
Choosing the parameters $a = q /B,\ b = q/ A,\ c
= q/ D, \ d = Cq\lambda_+$ and $e = Cq\lambda_-$ and
renormalizing we again obtain the equation (2.1) with a
solution
$$
\leqalignno{
X_n^{(6)} \ = &\ (-1)^n \big({q\over ABD}\big)^n
q^{-n(n-1)/2} \big( {ABD\over q} \lambda_+ , {ABD\over q}
\lambda_- \big)_n  & (2.19) \cr
& \quad \times \tphit\ \left( \matrix{ q^{-n+1}/B ,
 q^{-n+1}/A,  q^{-n+1}/D \cr 
C\lambda_+ q^{-n+1} , C\lambda_- q^{-n+1}\cr} ; Cq^n
\right)\cr}
$$
and three similar solutions obtained by parameter
interchanges $C \leftrightarrow A$ or $C \leftrightarrow B$
or $C \leftrightarrow D$. However, the solution
$X_n^{(6)}$ and its $C \leftrightarrow A,\ C \leftrightarrow
B$ or $C \leftrightarrow D$ interchanges do not give new
solutions.  They are related to the solutions $X_n^{(2)}, 
X_n^{(3)}, X_n^{(4)}$ and $X_n^{(5)}$ by the transformation
formula (2.2b).  For example, with the help of this formula
we find that
$$
X_n^{(6)} \ = \ {(Cq /A, Cq /D, Cq /B
)_\infty \over (C, C\lambda_+ q, C\lambda_- q)_\infty}
X_n^{(4)}
$$
and thus $X_n^{(6)}$ is the same solution as $X_n^{(4)}$
except for a constant factor.

We next show that continuous dual $q$-Hahn polynomials ([3,
pp. 3, 28]) are obtained as a particular case of the above
solutions.  This is true for the solution $X_n^{(2)}$ if
$C = q$ or $X_n^{(6)}$ if $A,B$ or $D = q$.  That is why our
general case represents associated continuous dual $q$-Hahn
polynomials.  In order to show this for $X_n^{(2)} $ we
first apply the transformation formula [5, (III.34), p. 245]
connecting $\tphit$'s of type I and II
$$
\eqalign{
\tphit\ \left( \matrix{a,b,c\cr d,e\cr}; {de\over abc}
\right) \ = & {({e/b}, {e/c})_\infty\over (e,
{e/bc})_\infty} \tphit \left(\matrix{ {d/a}, b,
c\cr d, {bcq/e}\cr} ; q\right)\cr
& + {({d/a}, b,c,{de/bc} )_\infty\over
(d,e,{bc/e}, {de/abc} )_\infty} \tphit 
\left(\matrix{{e/b}, {e/c}, {de/abc}\cr
{de/bc}, {eq/bc}\cr} ; q\right).\cr}
$$
We have from (2.11)
$$
\leqalignno{
\qquad \qquad X_n^{(2)} \ = & \ {(A,B)_n\over (AB)^n} \left[{(AC\lambda_-, AC
\lambda_+ )_\infty\over ({Aq/D}, {C/B} )_\infty}\
\tphit \ \left( \matrix{ q^{-n+1}/C, AB\lambda_+, AB\lambda_-\cr
{Aq/C}, {Bq/C}\cr}; q \right) \right. & (2.20) \cr
& +\left. {(q^{-n+1}/C, AB\lambda_+, AB\lambda_-, {Aq/B} )_\infty\over
({Aq/C}, {Aq/D}, {B/C}, q^{-n+1}/B )_\infty}
\ \tphit\ \left(\matrix{
 q^{-n+1}/B, AC\lambda_-, AC\lambda_+\cr
{Aq/B}, {Cq/B}\cr}; q \right)\right]\ .\cr}
$$
When we write $C = q$ in (2.20), the right side becomes
$$
{(A,B)_n\over (AB)^n} {(Aq\lambda_- , Aq\lambda_+ )_\infty\over
({Aq\over D}, {q\over B})_\infty} \ \tphit \ \left( 
\matrix{q^{-n}, AB\lambda_+, AB\lambda_-\cr A,B\cr} ; q\right).
\leqno(2.21)
$$
In order to compare this with continuous dual $q$-Hahn polynomials we
apply a transformation which will ultimately change the 
interval of orthogonality from $[- 1 /\alpha, 1 /\alpha ]$ when
$\ \alpha = \frac12 \sqrt{ABCD /q}$ is real, to $[-1,1]$.  Take
$x = \alpha z = \cos \theta , \ u = e^{i\theta}$, which
means $u = 2\alpha \lambda_+ = 1 /2\alpha \lambda_-$.
Thus, omitting constant factors we can write (2.21) as
$$
{(A,B)_n\over (AB)^n} \ \tphit \ \left( 
\matrix{q^{-n}, ABu /2\alpha , AB /2\alpha u\cr 
A,B\cr} ; q\right)  \leqno(2.22)
$$
which except for a normalization factor of $(2\alpha )^{-n}$, is the
same as the continuous dual $q$-Hahn polynomials of Askey
and Wilson [3].

We now proceed to show that $X_n^{(1),-}$ is a minimal
solution of (2.1) for a particular branch in the complex plane.
This is done by evaluating the large-$n$ asymptotics of 
$X_n^{(1),\pm}$ and the equation (2.1).  Applying the
transformation (2.2b) to (2.8) i.e. to $X_n^{(1),\pm}$ we
obtain for $n \to \infty$
$$
\leqalignno{
X_n^{(1),\pm} &\ = \ {(A,B,C,D)_n\over (BCD\lambda_\pm,
ABC\lambda_\pm )_n} (\lambda_\pm )^n {(Bq^n,ACD\lambda_\pm
q^n, ABCD\lambda_\pm^2 )_\infty\over (BCD\lambda_\pm q^n,
ABC\lambda_\pm q^n, AD\lambda_\pm )_\infty}\cr
& \qquad \times \tphit \ \left( \matrix{CD\lambda_\pm ,
AC\lambda_\pm , AD\lambda_\pm\cr ACD\lambda_\pm q^n,
ABCD\lambda_\pm^2\cr} ; Bq^n \right) & (2.23)\cr
& \qquad \approx \hbox{const } \ (\lambda_\pm )^n.\cr}
$$
On the other hand, asymptotics of the second order
difference equation (2.1) is given by
$$
X_{n+1} - zX_n + {q\over ABCD} X_{n-1} \ = \ 0, \leqno(2.24)
$$
from which we have for large $n$
$$
X_n \ \approx \hbox{const } \ (\lambda_\pm )^n. \leqno(2.25)
$$
Choosing the square root branch for which $| \lambda_-/
\lambda_+ | < 1$, the minimal (or the sub-dominant)
solution in terms of the variable $x = \alpha z$ is given by
$$
\leqalignno{
\qquad X^{(\min )}_n (x) &\ = \ X_n^{(1),-} & (2.26) \cr
&\ = \ (\lambda_- )^n {(A,B,C,D)_n\over (BCD\lambda_- ,
ABC\lambda_- )_n} \ \tphit\ \left( \matrix{Bq^n, Cq^n,
BC\lambda_-\cr BCD\lambda_- q^n, ABC\lambda_- q^n\cr} ;
AD\lambda_- \right)\cr}
$$
valid for $z\alpha = x \in {\Bbb C} \setminus [-1,1],\ \alpha =
\frac12 \sqrt{ABCD\over q}$, $\lambda_\pm = {1\over 2\alpha}
(x\pm \sqrt{x^2-1})$.  Summarizing what we have done so far, we have

\proclaim Theorem 1.  The functions 
$X_n^{(1),\pm}$ , $X_n^{(k)},\ k =
2,3,4,5$ of {\rm (2.8), (2.11), (2.13), (2.14)} and {\rm (2.15)}
respectively are solutions to the recurrence relation {\rm (2.1)} for
associated continuous dual $q$-Hahn polynomials. These solutions are 
pairwise linearly
independent. The minimal
solution of {\rm (2.1)} is given, up to a multiplicative factor, 
by {\rm (2.26)} with the square root
branch chosen so that $(\lambda_- 
/\lambda_+)^n \to 0$ as $n \to \infty$ with $z \alpha
\ = \ x \in {\Bbb C} \setminus [-1,1],\ \alpha = \frac12
\sqrt{ABCD /q}$.
\par
\vfill\eject

\noindent{\bf 3. The continued fraction and measure.}
\bigskip

The infinite continued fraction associated with (2.1) is
$$
CF(z) \ = \ z - a_0 - {b_1^2\over z-a_1}{{}\atop -}{b_2^2\over z-a_2}{{}\atop - \cdots},\
b_n^2 \not= 0,\ n > 0. \leqno(3.1)
$$
Pincherle's theorem [16] connects the minimal solution of (2.1)
with the continued fraction (3.1) by the formula
$$
{1\over CF(z)} \ = \ {X_0^{(\min )} (x)\over b_0^2 X_{-1}^{(\min )}
(x)},\ z\alpha \ = \ x,\ \alpha \ = \ \frac12
\sqrt{ABCD /q}. \leqno(3.2)
$$
Therefore from Theorem 1, we obtain the continued fraction
representation
$$
{1\over CF(z)} \ = \ {ABCD\lambda_-\over q(1-
{BCD\lambda_-\over q})(1-{ABC\lambda_-\over q})}\ 
{\tphit \ \left(\matrix{BC\lambda_-, B,C\cr BCD\lambda_-,
ABC\lambda_- \cr}; AD\lambda_- \right)\over 
\tphit \ \left(\matrix{BC\lambda_-, B /q ,C /q \cr
BCD\lambda_- /q , ABC\lambda_-  /q \cr}; 
AD\lambda_- \right)}. \leqno(3.3)
$$
We can also write (3.3) as
$$
{1\over CF(z)} \ = \ {1\over \lambda_+ (1-{1\over
A\lambda_+} )(1- {1\over D\lambda_+})} \
{\tphit \ \left(\matrix{BC\lambda_-, B,C\cr  
q /A\lambda_+  , q /D\lambda_+\cr}; AD\lambda_-
\right) \over 
\tphit \ \left(\matrix{BC\lambda_-, B /q ,C /q \cr
1 /A\lambda_+  , 1 /D\lambda_+\cr}; AD\lambda_-
\right)}. \leqno(3.3a)
$$
(3.3) or the alternative form (3.3a) are valid for
$$
z\alpha = x \in C \setminus [-1,1],\ | \lambda_-
/\lambda_+  | < 1,\ \alpha = \frac12 \sqrt{ABCD /q}.
$$
In the particular case $C = q$ (the case of continuous dual
$q$-Hahn polynomials), (3.3) reduces to
$$
{1\over CF(z)} \ = \ {ABD \lambda_-\over (1-BD \lambda_- )
(1-AB\lambda_-)} \ \tphit \ \left( \matrix{Bq\lambda_-,
B,q\cr DBq\lambda_-, ABq\lambda_-\cr} ; AD\lambda_- \right)
\leqno(3.4)
$$
which can also be written with the help of (2.2b) in the
form
$$
{1\over CF(z)} \ = \ {ABD \lambda_- (q,ABD\lambda_-,
ABDq\lambda_-^2 )_\infty\over (BD\lambda_-, AB\lambda_-, 
AD\lambda_-)_\infty} \ \tphit \ \left(\matrix{BD\lambda_-,
AB\lambda_- , AD\lambda_-\cr ABD\lambda_-, ABDq\lambda^2_-\cr}
;q\right). \leqno(3.5)
$$
with explicit pole terms given by the zeros of the 
denominator $(BD\lambda_-, AB\lambda_-, AD\lambda_- )_\infty$.
These pole singularities and their residues determine the discrete
component of the spectral measure of orthogonality for continuous
dual $q$-Hahn polynomials.
\bigskip

We now determine the absolutely continuous part of the spectrum for
the general associated case.

If $x =z\alpha \in (-1,1)$, then (2.1) has linearly independent solutions
given by the boundary values of (2.26) as $x$ approaches $(-1,1)$ from
above and below. With now $\lambda_\pm = 
{1\over 2\alpha} [x\pm i\sqrt{1-x^2}], \ x = z\alpha$, we have the
large $n$ asymptotics
$$
X_n^{(\min )} (x+i0) \ \approx \ (\lambda_- )^n {(A,B,C,D)_n (ABCD
\lambda^2_- )_\infty\over (BCD\lambda_-, ABC\lambda_-, AD\lambda_-)
_\infty}, \leqno(3.6)
$$
and
$$
X_n^{(\min )} (x-i0) \ \approx \ (\lambda_+ )^n {(A,B,C,D)_n (ABCD
\lambda^2_+ )_\infty\over (BCD\lambda_+, ABC\lambda_+, AD\lambda_+)
_\infty}. \leqno(3.7)
$$
Since the minimal solution changes as we cross the line segment
$z\alpha =x \in (-1,1)$, we have the representation
$$
{1\over CF(z)} \ = \ \int_\bR {\omega(t)dt\over z-t/\alpha }
+ \ \hbox{possible pole terms}, \leqno(3.8)
$$
Also $ \omega (x) ,\ x \in (-1,1)$ can
be obtained by using the formula [21]
$$
\leqalignno{
\omega (x) \ = \ & {1\over 2\pi i \alpha} \ {W (X_{-1}^{(\min)} (x+i0),
X_{-1}^{(\min)} (x-i0))\over b_0^2  X_{-1}^{(\min)} (x+i0)X_{-1}^{(\min)}(x-i0)} & (3.9) \cr
\ = \ & {1\over 2\pi i \alpha} \lim\limits_{n\to \infty}
\ {W (X_{n}^{(\min)} (x+i0),
X_{n}^{(\min)} (x-i0))\over b_1^2 b_2^2 \ldots b_n^2 b_0^4 
 X_{-1}^{(\min)} (x+i0)X_{-1}^{(\min)} (x-i0)},\cr}
$$
where
$$
W(X_n,Y_n) \ = \ X_n Y_{n+1} -X_{n+1}Y_n.
$$
Using (2.1), (2.26), (3.6) and (3.7) and simplifying we have
$$
\leqalignno{
& \qquad \omega (x) \ = \ {1\over 2\pi \sqrt{1-x^2}} {(A,B,C,D)_\infty 
({1\over u^2},u^2)_\infty\over ({2\alpha\over Au},
{2\alpha\over A}u, {2\alpha \over Du}, {2\alpha\over D}u,
{2\alpha q\over BCu}, {2\alpha q\over BC}u)_\infty } & (3.10) \cr
& \quad \times \left[ \tphit\ \left(\matrix{BC\lambda_- , {B\over
q}, {C\over q}\cr {BCD\lambda_-\over q}, {ABC\lambda_-\over
q}\cr}; AD\lambda_-\right)\  
\tphit\ \left(\matrix{BC\lambda_+ , {B\over
q}, {C\over q}\cr {BCD\lambda_+\over q}, {ABC\lambda_+\over
q}\cr}; AD\lambda_+\right)  \right]^{-1}. \cr}
$$
In the particular case $C = q$ this reduces to
$$
\omega (x) \ = \ {1\over 2\pi \sqrt{1-x^2}} \>
{(A,B,q,D)_\infty ({1\over u^2},u^2)_\infty\over
(\sqrt{BD\over A}/u, \sqrt{BD\over A} u , \sqrt{AB\over
D}/u, \sqrt{AB\over D}u, \sqrt{AD\over B}/u, \sqrt{AD\over
B}u)_\infty }. \leqno(3.11)
$$

Taking the appropriate value of the parameters from (2.22) we
find that this weight function is the same as the one obtained by
Askey and Wilson ([3, p. 11, Theorem 2.2 in the special case
$d= 0$]) for continuous dual $q$-Hahn polynomials.
Summarizing the above we have

\proclaim Theorem 2.  The associated continuous dual
$q$-Hahn polynomials $P_n (x /\alpha)
$ given by {\rm (4.12)} of the next
section are orthogonal with respect to a measure with
absolutely continuous component given by the weight function
{\rm (3.10)} on $(-1,1)$.  In the particular case $C = q$ this
absolutely continuous component reduces to {\rm (3.11)} and the
discrete spectrum is given by the zeros of $(BD\lambda_-,
AB\lambda_-, AD\lambda_- )_\infty$ where $\lambda_- =
(x \mp \sqrt{x^2-1})/2\alpha$ for $x>1$ and $x<-1$ respectively.
\par

The $q\to 1$ limit of continuous dual $q$-Hahn polynomials yields
the case of continuous dual Hahn polynomials. The corresponding
results of this Section for associated continuous dual Hahn
polynomials with $q=1$ are given by Ismail, Letessier and Valent [13].
\vfill\eject

\noindent{\bf 4. Generating function.}
\bigskip

The associated continuous dual $q$-Hahn polynomials $P_n
(z;A,B,C,D)$ satisfy the second order difference equation
(2.1) i.e. the equation
$$
\leqalignno{
& X_{n+1} - \left[z- \left({1\over A} + {1\over B} + {1\over C} +
{1\over D} \right) q^n + (1+q) q^{2n-1} \right] X_n  & (4.1) \cr
& \qquad + {1\over 4\alpha^2} (1-Aq^{n-1} )(1-Bq^{n-1})
(1-Cq^{n-1})(1-Dq^{n-1} ) X_{n-1} \ = \ 0.\cr}
$$
A renormalized form of (4.1) is, with $x = \alpha z,\ \alpha
\ = \ \frac12 \sqrt{ABCD /q}$,
$$
\displaylines{
(4.2)\quad (1-Aq^n)(1-Dq^n)\zeta_{n+1} - \left[2x-2\alpha \left({1\over A} +
{1\over B} +{1\over C} + {1\over D}\right) q^n + 2\alpha (1+q)
q^{2n-1} \right] \zeta_n\hfill\cr
\hfill {} + (1-Bq^{n-1})(1-Cq^{n-1} )\zeta_{n-1} \
= \ 0. \cr}
$$
This is satisfied by the polynomials
$$
\zeta_n (x; A,B,C,D) \ = \ {(2\alpha)^n P_n (z;A,B,C,D)\over
(A)_n (D)_n}.
$$
Let the generating function of the polynomials $\zeta_n$ be
$$
\leqalignno{
G(x,t) \ = & \ \sum_{n=0}^\infty {(2\alpha )^n P_n (z; A,B,C,D)
\over (A)_n (D)_n} t^n & (4.3) \cr
\ = & \ \sum_{n=0}^\infty \zeta_n (x;A,B,C,D)t^n.\cr}
$$
An explicit form for the generating function may be obtained by
employing the procedure given in [14].  First multiply (4.2)
by $t^n$ and sum up the resulting equations from $n = 0$ through
$n = \infty$.  Using the initial conditions for the polynomial
solutions of the first kind, we obtain
$$
\leqalignno{
& (1-ut)(1-{t\over u}) G(x,t) -(1-{2\alpha tq\over AD} )({A\over q}
+ {D\over q} - {2\alpha t\over B} - {2\alpha t\over C}) G(x,tq) &
(4.4) \cr
& \qquad + (1-{2\alpha tq\over AD} ) (1-{2\alpha tq^2\over AD} ) 
{AD\over q^2} G(x,tq^2) \ = \ (1-{A\over q})(1-{D\over q}).\cr}
$$
In (4.4) we put
$$
G(x,t) \ = \ {\left({2\alpha tq\over AD}\right)_\infty\over
(tu)_\infty} f(t) \leqno(4.5)
$$
to obtain
$$
\leqalignno{
& (1-{t\over u} ) f(t) - ({A\over q} + {D\over q} - {2\alpha t\over B}
- {2\alpha t\over C} ) f(tq) & (4.6)\cr 
& \qquad + {AD\over q^2} (1-tuq) f(tq^2) \ = \ (1-{A\over q})
(1-{D\over q}) {(tuq)_\infty\over \left({2\alpha tq\over AD} 
\right)_\infty} .\cr}
$$
In the left side of (4.6) we write $f(t) = \sum_{n=0}^\infty f_nt^n$
and on the right side we use $q$-binomial theorem to replace
$ (tuq)_\infty /\left({2\alpha tq\over AD} \right)_\infty $
by
$$
\sum_{n=0}^\infty \ {\left({ADu\over 2\alpha}\right)_n\over
(q)_n} \left({2\alpha tq\over AD} \right)^n .
$$
If we now equate coefficients of $t^n$ on both sides, we
obtain the first order difference equation
$$
\displaylines{
(4.7)\quad (1-Aq^{n-1})(1-Dq^{n-1}) f_n \ = \ {1\over u} (1-
{2\alpha u\over B} q^{n-1} )(1-{2\alpha u\over C} q^{n-1})
f_{n-1} \hfill\cr 
\hfill {} + \left({2\alpha q\over AD} \right)^n (1-{A\over
q}) (1- {D\over q} ) {\left({ADu\over 2\alpha}\right)_n\over
(q)_n}.  \cr}
$$
Rewriting (4.7) as
$$
\leqalignno{
u^n {(A)_n (D)_n\over \left({2\alpha u\over B}\right)_n
\left({2\alpha u\over C} \right)_n } f_n \ = &\  u^{n-1} 
{(A)_{n-1} (D)_{n-1}\over \left({2\alpha u\over B}
\right)_{n-1} \left({2\alpha u\over C} \right)_{n-1}} f_{n-1} & (4.8)\cr
& + \left({2\alpha q\over AD}\right)^n u^n {({A\over q})_n
({D\over q})_n ({ADu\over 2\alpha})_n\over ({2\alpha u\over
B})_n ({2\alpha u\over C} )_n (q)_n}, \cr}
$$
the general solution of (4.7) is
$$
f_n \ = \ {({2\alpha u\over B})_n ({2\alpha u\over
C})_n\over u^n (A)_n (D)_n} \left[ E + \sum_{j=0}^n 
{({A\over q})_j ({D\over q})_j ({ADu\over 2\alpha})_j\over 
({2\alpha u\over B})_j ({2\alpha u\over C})_j (q)_j} \left(
{2\alpha q\over AD} \right)^j u^j \right], \leqno(4.9)
$$
where $E$ is a constant which by the boundary conditions may
be taken as 0.  Consequently we have the generating 
function 
$$
\leqalignno{
& \sum_{n=0}^\infty {(2\alpha)^n P_n (z;A,B,C,D)\over (A)_n (D)_n}
t^n  & (4.10) \cr
&\ = \ {({2\alpha t q\over AD} )_\infty\over (tu)_\infty}
\sum_{n=0}^\infty \sum_{j=0}^n {({2\alpha u\over B})_n ({2\alpha
u\over C})_n\over (A)_n (D)_n} {({A\over q})_j ({D\over q})_j 
({ADu\over 2\alpha})_j\over ({2\alpha u\over B})_j ({2\alpha u
\over C})_j (q)_j} \left({2\alpha q\over AD}\right)^j
u^{-n+j} t^n,\cr}
$$
where $z = x /\alpha$ and $x = \cos \theta,\ u = e^{i\theta}$.
If we interchange $C \leftrightarrow D$ in (4.10) and write
$C = q$ we obtain, with $a=2\alpha /B, b=2\alpha /D, c=2\alpha /A$,
$$
\leqalignno{
\sum_{n=0}^\infty &\  {(2\alpha)^n P_n (z; A,B,C)\over
(A)_n (q)_n} t^n  \cr
&\ = \ {(ct)_\infty\over (tu)_\infty} \sum_{n=0}^\infty
{(au)_n (bu)_n\over (q)_n (ab)_n} \left({t\over u}\right)^n &(4.11)
\cr
&\ = \ {(ct)_\infty\over (tu)_\infty} \ {}_2\phi_1 \left(\matrix{
au,bu\cr ab\cr} ; {t\over u} \right).\cr}
$$
This gives the
generating function result for the continuous dual $q$-Hahn
polynomials (see [13,  (3.3.7), p. 55]).

By comparing coefficients of $t^n$ on the left and right 
sides of (4.10) we obtain an explicit expression for our
monic associated continuous dual $q$-Hahn polynomials.
$$
\leqalignno{
P_n (z) \ = &\ P_n (z;A,B,C,D)\cr
\ = &\ \left({u\over 2\alpha}\right)^n {(A,D, {2\alpha q\over ADu})_n\over (q)_n}
\sum_{\ell = 0}^n \left\{ {(q^{-n}, {2\alpha u\over B}, {2\alpha u
\over C} )_\ell\over ({ADu\over 2\alpha} q^{-n}, A,D)_\ell}
\left({AD\over 2\alpha u}\right)^\ell \right. & (4.12)  \cr
& \qquad \left. \times \sum_{j=0}^\ell {({A\over q} , {D\over q} , 
{ADu\over 2\alpha})_j\over (q, {2\alpha u\over B}, {2\alpha u
\over C})_j} \left( {2\alpha uq\over AD}\right)^j \right\},\cr
\cr
\alpha \ = &\ \frac12 \sqrt{ABCD /q},\quad x \ = \ \frac12
(u+1/u), \quad z \ = \ x/\alpha \ .\cr}
$$
Note that $P_n (z)$ is of course symmetric in the parameters
$A,B,C,D$. The symmetry under the interchanges $A \leftrightarrow D$
or $B \leftrightarrow C$ is obvious from (4.12). However, the symmetry 
under the interchanges $A \leftrightarrow B,
A \leftrightarrow C, D \leftrightarrow B$ or $D \leftrightarrow C$
is hidden. $P_n(z)$ is also symmetric under the interchange
$u \leftrightarrow u^{-1}$. Again this is not apparent from (4.12).
Applying any one of these hidden symmetry interchanges to (4.12)
gives us a type of transformation formula.

A different expression for $P_n(z)$ is obtained in Section 7. It is
derived from the associated Askey-Wilson polynomial formula of
Ismail and Rahman [15]. In order to contrast (7.3) with (4.12) we
repeat it here as
$$
\leqalignno{
P_n(z; A,B,C,D) \ = \ & {(B,C)_n\over (BC)^n} \left\{
\sum_{k=0}^n {(q^{-n}, \sqrt{{BCq\over AD}} u,
\sqrt{{BCq\over AD}} {1\over u})_k \ q^k\over (q,B,C)_k}\right.
& (4.13)\cr
& \left. \times \sum_{j=0}^{n-k} {({A\over q}, {D\over q}, q^{k+1},
q^{k-n})_j\over (q,Cq^k , Bq^k, q^{-n})_j} \left({BCq\over
AD}\right)^j \right\}.\cr}
$$
This formula also does not reveal the full symmetry with respect to
$A,B,C,D$. However, it does make explicit the $u\leftrightarrow u^{-1}$ 
symmetry.
\vfill\eject
\noindent{\bf 5. Four limiting cases.}
\bigskip

We now take successive limits $D \to \infty,\ C \to \infty,
\ B \to \infty,\ A \to \infty$.
\bigskip

\noindent{\bf 5.1 Associated big $q$-Laguerre ($D\to \infty$).}
\bigskip

The recurrence relation (2.1) becomes
$$
\leqalignno{
L_{n+1} - & \left[z - \left({1\over A} + {1\over B} + {1\over C} \right)q^n + (1+q)
q^{2n-1}\right] L_n & (5.1) \cr
& \qquad - {q^n\over ABC} (1-Aq^{n-1})(1-Bq^{n-1}) (1-Cq^{n-1})
L_{n-1} \ = \ 0.\cr}
$$

If we write $A = aq,\ B = bq,\ C = -abq /s$ and
$z = t /abq$, and renormalize with a factor
$(abq)^n$, (5.1) gives the equation (4.3) of Ismail and Libis [14]
for associated big $q$-Laguerre polynomials in monic form.

The solutions of (5.1) are obtained as $D \to \infty$ limits of
solutions of (2.1).  We have $\lambda_+ \approx z$ and
$\lambda_- \approx q/ABCDz$ and consequently from (2.8)
$$
\leqalignno{
L_n^{(1)} (z;A,B,C) \ = \ & \lim_{D\to\infty} X_n^{(1),-} (z;
A,B,C,D) & (5.2) \cr
\ = \ & (-1)^n {q^{n(n+1)/2}\over (ABCz)^n} {(A,B,C)_n
\over ({q\over Az})_n} \ {}_2\phi_1\ \left(
\matrix{Bq^n, Cq^n\cr q^{n+1}/ Az\cr} ; {q\over BCz} \right),\cr}
$$
and from (2.12)
$$
\leqalignno{
L_n^{(2)} (z;A,B,C) \ = \ & {(q/B)_\infty \over(Aq/B)_\infty }
\lim_{D\to \infty} X_n^{(2)} (z; A,B,C,D) & (5.3)\cr
\ = \ & (-1)^n q^{n(n-1)/2}(A)_nA^{-n} \ {}_2\phi_2
\ \left(\matrix{Aq^n, {q\over BCz}\cr {Aq\over B}, {Aq\over C}\cr};
Azq^{-n+1}\right).\cr}
$$

Interchanging $A \leftrightarrow B$ and $A \leftrightarrow C$ in (5.3) yields
the limits of solutions $X_n^{(3)},\ X_n^{(4)}$.  We have
$$
L_n^{(3)} (z; A,B,C) \ = \ (-1)^n q^{n(n-1) /2} (B)_n B^{-n}
{}_2\phi_2 \ \left(\matrix{Bq^n, {q\over ACz}\cr {Bq\over C},
{Bq\over A}\cr}; Bzq^{-n+1} \right), \leqno(5.4)
$$
and
$$
L_n^{(4)} (z; A,B,C) \ = \ (-1)^n q^{n(n-1) /2} (C)_nC^{-n}
{}_2\phi_2 \ \left(\matrix{Cq^n, {q\over ABz}\cr {Cq\over A},
{Cq\over B}\cr}; Czq^{-n+1} \right). \leqno(5.5)
$$
Limits of $X_n^{(6)}$ given by (2.19) and its parameter interchanges
will be in terms of ${}_2\phi_1$ and these three limits will
simply be transforms of $L_n^{(2)}, L_n^{(3)}$ and $L_n^{(4)}$.
We write below one of these limits 
$$
\leqalignno{
L_n^{(5)} (z;A,B,C) \ = &\ \lim\limits_{D\to \infty} X_n^{(6)} 
(z; A,B,C,D) \cr
\ = &\ z^n (1 /Cz)_n \ {}_2\phi_1 \ \left(
\matrix{ q^{-n+1}/A, q^{-n+1}/B\cr Cz q^{-n+1}\cr} ; Cq^n\right).
& (5.6)\cr}
$$

Next we take the $D \to \infty$ limit of the explicit expression
(4.12) for the polynomial solution to get
$$
\leqalignno{
P_n^{(1)} (z; A,B,C) \ = \ & \lim_{D\to \infty} P_n (z;A,D,C,B) 
& (5.7)\cr
\ = \ & z^n {(A,B, {q\over ABz})_n\over (q)_n} \sum_{\ell=0}^n
\left\{ {(q^{-n}, {ABCz\over q})_\ell (-1)^\ell q^{\ell (\ell -1)/2}
\over (ABzq^{-n}, A,B)_\ell} \left({AB\over C}\right)^\ell\right.\cr
& \quad \left. \times \sum_{j=0}^\ell {({A\over q}, {B\over q}, ABz)_j\over
({ABCz\over q}, q)_j} (-1)^j q^{-j(j-1)/2} \left({Cq\over AB}
\right)^j \right\}.\cr}
$$
Note that before taking limit $D \to \infty$ of (4.12) we have
made the parameter interchange $B \leftrightarrow D$ in 
(4.12).

The minimal solution of (5.1) is $L_n^{(1)}$ and we therefore
have the continued fraction representation
$$
\leqalignno{
{1\over z-a_0}{{}\atop -} {b_1^2\over z-a_1}{{}\atop -} {b_2^2\over z-a_2}{{}\atop- \cdots} & \ = \
{{}_2\phi_1 \left(\matrix{B,C\cr q/Az\cr}; {q\over BCz}\right)\over
z(1-1/Az ) {}_2\phi_1 \left(\matrix{B/q, C/q\cr 1/Az\cr} ; {q\over BCz}
\right)}, & (5.8) \cr}$$
$$a_n \ = \ ({1\over A} + {1\over B} + {1\over C} ) q^n - (1+q) q^{2n-1}
$$
$$
b_n^2  \ = \ - {q^n\over ABC} (1-Aq^{n-1}) (1-Bq^{n-1} )(1-Cq^{n-1}).
$$
The associated orthogonality is discrete and only explicit in the 
case $A,B$ or $C = q$.
\vfill\eject

\noindent{\bf 5.2.  Associated Wall $(C,D \to \infty)$.}

With $C\to \infty$ the recurrence relation (5.1) changes to
$$
W_{n+1} - \left[z-\left({1\over A} + {1\over B} \right)q^n + (1+q) q^{2n-1}\right]W_n + 
{q^{2n-1}\over AB} (1-Aq^{n-1})(1-Bq^{n-1}) W_{n-1} \ = \ 0.
\leqno(5.9)
$$
Solutions are
$$
\leqalignno{
W_n^{(1)} (z;A,B) \ = \ & \lim\limits_{C\to \infty} L_n^{(1)} 
(z;A,B,C) & (5.10) \cr
\ = \ & \left({q\over ABz} \right)^n q^{n(n-1)} {(A,B)_n\over
(q /Az)_n} \ {}_1\phi_1 \left(\matrix{Bq^n\cr q^{n+1} /Az
\cr}; {q^{n+1}\over Bz} \right),\cr}
$$
$$
\leqalignno{
W_n^{(2)} (z;A,B) \ = \ & \lim\limits_{C\to \infty} L_n^{(2)} 
(z;A,B,C) & (5.11) \cr
\ = \ & (-1)^n  q^{n(n-1)/2} (A)_n A^{-n}
\ {}_1\phi_1 \left(\matrix{Aq^n\cr Aq /B \cr}; 
Azq^{-n+1} \right),\cr}
$$
and
$$
\leqalignno{
W_n^{(3)} (z;A,B) \ = \ & \lim\limits_{C\to \infty} L_n^{(3)} 
(z;A,B,C) & (5.12) \cr
\ = \ & (-1)^n  q^{n(n-1)/2} (B)_n B^{-n}
\ {}_1\phi_1 \left(\matrix{Bq^n\cr Bq /A \cr}; 
Bzq^{-n+1} \right),\cr}
$$
which is just $A \leftrightarrow B$ interchange of (5.11).

Also from (5.6)
$$
\leqalignno{
W_n^{(4)} (z;A,B) \ = \ & \lim\limits_{C\to \infty} L_n^{(5)}(z;A,B,C) & (5.13) \cr
\ = \ & z^n {}_2\phi_0 \left(\matrix{q^{-n+1}/A, q^{-n+1}/B \cr
-\cr} ; {q^{2n-1}\over z} \right).\cr}
$$
The series representing the above ${}_2\phi_0$ converges only when it
terminates. the relevant terminating cases are when $A=q$ or $B=q$.

A limit $C \to \infty$ of (5.7) gives explicit expression for
associated Wall polynomials
$$
\leqalignno{
P_n^{(2)} (z;A,B) \ = \ & \lim\limits_{C\to \infty} P_n^{(1)} (z;A,B,C) & (5.14)\cr
\ = \ & z^n {(q/ABz, A,B)_n\over (q)_n} \sum_{\ell = 0}^n
\left\{ {(q^{-n})_\ell q^{\ell (\ell -1)}\over (q^{-n} ABz, A,B)_\ell}
\left({A^2B^2z\over q}\right)^\ell\right.\cr
& \qquad \left. \times \sum_{j=0}^\ell {(A/q, B/q, ABz)_j\over (q)_j} 
\left({q\over AB}\right)^{2j} z^{-j} q^{-j(j-1)} \right\} .\cr}
$$
The minimal solution of (5.9) is given by (5.10) and
therefore we have the associated continued fraction representation
$$
{1\over z-a_0}{{}\atop -} {b_1^2\over z-a_1}{{}\atop -} 
{b_2^2\over z-a_2}{{}\atop- \cdots}  \ = \
{{}_1\phi_1 \left(\matrix{B\cr q/Az\cr}; {q\over Bz}\right)\over
z(1-1/Az ) {}_1\phi_1 \left(\matrix{B/q \cr 1/Az\cr} ; {1\over Bz}
\right)}, \leqno(5.15) 
$$
where
$$
\eqalign{
a_n & \ = \ ({1\over A} + {1\over B} ) q^n - (1+q) q^{2n-1}\cr
b_n^2 & \ = \  {q^{2n-1} \over AB} (1-Aq^{n-1}) (1-Bq^{n-1} ).\cr}
$$
When $A$ or $B = q$, (5.15) has associated with it a discrete 
orthogonality for $P_n^{(2)} (z)$.  We note here that our
associated Wall polynomials given by (5.9) reduce to the
Wall polynomials when $B = q$.  In fact, if we make the
substitutions $A = aq,\ B = q, \ z = x /aq$ in (5.9)
and renormalize, the equation can be written in the 
form (see [4])
$$
q^n (1-aq^{n+1})p_{n+1}(x) - \left[q^n (1-aq^{n+1} ) + aq^n (1-q^n)-x\right]
p_n(x) + aq^n (1-q^n) p_{n-1}(x) = 0
$$
with 
$$
p_n(x)=p_n (x;a) = (-1)^n q^{-n(n-1)/2} {(aq)^n\over (aq)_n} W^{(2)}_n 
({x\over aq}; aq,q)
$$
where we use solution $W_n^{(2)}$ from (5.11). From (5.28) we obtain the 
standard expression for Wall polynomials (see [17, (3.20.1), p. 83] and [4, p. 198])
$$
p_n (x;a) \ = \ {}_2\phi_1 \left(\matrix{q^{-n},0\cr aq\cr}; qx\right).
$$
\bigskip

\noindent{\bf 5.3. Limit Wall $(B,C,D \to \infty)$.}
\bigskip

The three-term recurrence now becomes
$$
U_{n+1} - \left[z - {q^n\over A} + (1+q) q^{2n-1} \right]U_n - {q^{3n-2}\over A}
(1-Aq^{n-1}) U_{n-1} \ = \ 0. \leqno(5.16)
$$
Using (5.10), (5.11) and (5.13), we have the solutions
$$
\leqalignno{
U_n^{(1)} (z;A) &\ = \ \lim\limits_{B\to \infty} W_n^{(1)} (z;A,B)
& (5.17) \cr
& \ = \ (-1)^n ({q\over Az})^n q^{3n(n-1)/2} {(A)_n\over ({q\over
Az})_n} {}_0\phi_1 \left(\matrix{ -\cr {q^{n+1}\over Az}\cr};
{q^{2n+1}\over z} \right),\cr
U_n^{(2)} (z;A) & \ = \ \lim\limits_{B\to \infty} W_n^{(2)} (z;A,B)
& (5.18) \cr
& \ = \ (-1)^n  q^{n(n-1)/2} {(A)_n\over A^n}
{}_1\phi_1 \left(\matrix{ Aq^n\cr 0\cr};
Azq^{-n+1} \right), \cr}
$$
and
$$
U_n^{(3)} (z;A) \ = \ \lim\limits_{B\to \infty} W_n^{(4)} (z;A,B)
\ = \ z^n  {}_2\phi_0 \left(\matrix{ {q^{-n+1}\over A},0\cr
-\cr}; {q^{2n-1}\over z} \right), \leqno(5.19)
$$
where $U_n^{(3)}$ converges when it terminates with say $A=q$.

A direct limit of (5.14) i.e. $P_n^{(2)} (z; A,B)$ as $B \to \infty$
leads to an indeterminate form.  However we can obtain the explicit form 
of the polynomials by applying the method of Section 4
abinitio to equation (5.16).  The result is
$$
\leqalignno{
P_n^{(3)} (z;A)
& \ = \ {q^{n^2}\over A^n} {(A)_n\over (q)_n} \sum_{\ell = 0}^n
\left\{ {(q^{-n})_\ell\over (A)_\ell} (-1)^\ell q^{-\ell (\ell -1)/2}
(Az)^{\ell} \right. & (5.20) \cr
&\left. \qquad \times \sum_{j=0}^\ell {(A /q)_j\over (q)_j}
q^{j(j-1)/2} (-Az)^{-j}\right\}.\cr}
$$
The minimal solution of (5.16) being (5.17) we have the related
continued fraction
$$
\leqalignno{
{1\over z-a_0}{{}\atop-} {b_1^2\over z-a_1}{{}\atop -} 
{b_2^2\over z-a_2}{{}\atop- \cdots}  \ = \ & {1\over z}
{{}_0\phi_1 \left(\matrix{-\cr q/Az\cr}; {q\over z}\right)\over
(1-1/Az ) {}_0\phi_1 \left(\matrix{- \cr 1/Az\cr} ; {1\over qz}
\right)}, & (5.21) \cr
 \ = \ & {1\over z}
{{}_1\phi_1 \left(\matrix{A\cr 0\cr}; {q\over Az}\right)\over
{}_1\phi_1 \left(\matrix{A/q \cr 0\cr} ; {1\over Az}
\right)}, \cr}
$$
where
$$
\eqalign{
a_n & \ = \ {q^n\over A}  - (1+q) q^{2n-1}\cr
b_n^2 & \ = \ - {q^{3n-2}\over A} (1-Aq^{n-1}) .\cr}
$$
The second expression on the right side of (5.21) comes from the
transformation
$$
{}_0\phi_1 \left( \matrix{-\cr c\cr} ; cz \right) \ = \
{1\over (c)_\infty} \ {}_1\phi_1 \left(\matrix{z\cr 0\cr} ; c\right)
$$
which can be derived from (5.30) by letting $a = 0$.  Note that
when $A \not= q$, the right side of (5.21) is a meromorphic
function of $z$.  When $A \to q$ the singularities 
coalesce at $z = 0$ to produce an essential singularity.
\bigskip

\noindent{\bf 5.4 A fourth limit $(A,B,C,D \to \infty)$.}
\bigskip

The three-term recurrence is now
$$
V_{n+1} - [z + (1+q) q^{2n-1} ] V_n + q^{4n-3} V_{n-1} \ = \ 0.
\leqno(5.22)
$$

Using (5.17) and (5.19) we have the solutions
$$
V_n^{(1)} (z)  \ = \ \lim_{A\to \infty} U_n^{(1)} (z;A) 
\ = \ q^{2n(n-1)} ({q\over z})^n {}_0\phi_1 \left(\matrix{
-\cr 0\cr}; {q^{2n+1}\over z} \right), \leqno (5.23)
$$
$$
V_n^{(2)} (z)  \ = \ \lim_{A\to \infty} U_n^{(3)} (z;A)
\ = \ z^n {}_2\phi_0 \left(\matrix{
0,0\cr -\cr}; {q^{2n-1}\over z} \right).\leqno (5.24)
$$
The solution $V_n^{(2)}$ is divergent and is thus only a formal
solution. The associated polynomials are given by
$$
\leqalignno{
P_n^{(4)} (z) \ = \ & \lim_{A\to \infty} P_n^{(3)} (z;A) & (5.25)\cr
\ = \ & {(-1)^n q^{n^2} q^{n(n-1)/2}\over (q)_n} \sum_{\ell = 0}^n 
\left\{ (q^{-n})_\ell q^{-\ell(\ell -1)} z^\ell \sum_{j=0}^\ell
{q^{j(j-1)}\over (q)_j} (qz)^{-j} \right\}.\cr}
$$
The minimal solution of (5.22), being given by (5.23), yields the associated 
continued fraction 
$$
{1\over z-a_0}{{}\atop-} {b_1^2\over z-a_1}{{}\atop-} {b_2^2\over z-a_2}{{}\atop- \cdots}
\ = \ {1\over z} {{}_0\phi_1 \left(\matrix{-\cr 0\cr}; q /z\right)\over
 {}_0\phi_1 \left(\matrix{-\cr 0\cr}; 1 /qz \right)} \leqno(5.26)
$$
where
$$
\leqalignno{
a_n &\ = \ -(1+q) q^{2n-1}\cr
b_n^2 &\ = \ q^{4n-3}.\cr}
$$
We can also write (5.26) more explicitly as
$$
{1\over z + (1+q)q^{-1}}{{}\atop -} {q\over z+(1+q)q}{{}\atop -} {q^5\over  
z+(1+q)q^3}{{}\atop -\cdots} 
 = {1\over z}{\sum_{n=0}^\infty {q^{n^2}\over (q)_n} z^{-n}\over
\sum_{n=0}^\infty {q^{n^2-2n}\over (q)_n} z^{-n}}. \leqno(5.27)
$$
Using (5.22), (5.23) and (5.26) we have the following

\proclaim Corollary.
If $\ 0 < q < 1$ and $n$ is an integer, then
$$
f_n (z) \ = \ {}_0\phi_1 \left(\matrix{-\cr 0\cr}; {q^{2n+1}\over
z} \right)
$$
has only real simple negative zeros which interlace those of $f_{n+1} (z)$.
\par
{\it Proof.} If $0<q<1$ then (5.26) is a completely convergent positive
definite $J$-fraction which can be represented as a Stieltjes
transform of a unique positive discrete measure [21]. This means that
(5.26) can have only simple pole sigularities on the real axis with 
positive residues. Hence we must have simple
real intertwining zeros for $f_0(z)$ and $f_{-1}(z)$. From the
series representation $f_n(z)=\sum_{k=0}^{\infty}{q^{k(k-1)}\over (q)_k}
(q^{2n+1}/z)^k$ 
we see that the zeros of $f_{-1}(z)$ and $f_0(z)$ must be negative.
This establishes the result for $n=-1$.
The proof for other values of $n$ is the same if one starts from
the continued fraction 
$$
{1\over z-a_{n+1}}{{}\atop -}{b_{n+2}^2\over z-a_{n+2}}{{}\atop -
\cdots}= {1\over z} {f_{n+1}(z)\over f_n(z)}.$$
\bigskip

\noindent{\it Note 1.} \quad There are similar corollaries
associated with the positive definite cases for the continued
fractions (5.8), (5.15) and (5.21).  These require special
parameter conditions. See also the remarks after (6.18), (6.28), (6.33) 
and (6.63).
\bigskip

\noindent{\it Note 2.} \quad The identities
$$
\leqalignno{
{}_2\phi_1 \left(\matrix {a,b\cr c\cr}; z\right) & \ = \ 
{(b,az)_\infty\over (c,z)_\infty} \ {}_2\phi_1 \left( 
\matrix{c/b,z\cr az\cr} ; b\right) & (5.28)\cr
{}_2\phi_1 \left(\matrix {a,b\cr c\cr}; z\right) & \ = \ 
{(az)_\infty\over (z)_\infty} \ {}_2\phi_2 \left( 
\matrix{a, c/b\cr c, az\cr} ; bz\right) & (5.29)\cr
{(z)_\infty\over (az)_\infty}\ {}_2\phi_1 \left(
\matrix {a,0\cr c\cr}; z\right) & \ = \ {}_1\phi_2 \left( 
\matrix{a\cr c, az\cr} ; cz\right) & (5.30)\cr
& \ = \  {1\over (c)_\infty} \ {}_1\phi_1 \left(\matrix{z\cr
az\cr} ; c\right)\cr}
$$
can be used to relate some of the above solutions.  For
(5.28) and (5.29) see [5, (III.1), (III.4), p. 241].
(5.30) follows from $b \to 0$ in (5.28)
and (5.29).
\vfill\eject

\noindent{\bf 6.  Additional Limits.}

There are other less obvious limiting cases which we can obtain from
(2.1) and its solutions. These may also be re-expressed as birth and
death processes. In the cases of Sections 6.4 and 6.5 there are  
process based seperately on the even and odd approximants with 
$z^2$ replaced by $z$ [7, Section 4]. We begin with
\medskip

\noindent{\bf  6.1  Associated Al-Salam-Chihara.}
\medskip

In (2.1) we put $D = \delta C$, multiply by $C$, replace $zC$ by
$z$ and renormalize and let $C \to 0$ to get
$$
Q_{n+1} - (z-(1+  \delta ^{-1}) q^n ) Q_n + {q\over AB\delta}
(1-Aq^{n-1} ) (1 - Bq^{n-1} ) Q_{n-1} = 0 
\leqno(6.1)
$$
with
$$
Q_n \ = \ \lim\limits_{C\to 0} C^n X_n \left( z /C; A,B,C,
C\delta \right). \leqno(6.2)
$$
This can be recognized as the recurrence relation for
associated Al-Salam-Chihara polynomials given in [2, (3.54) with $a=
(1+\gamma^2b)/\gamma$] if we make the replacements
$$(A,B,\delta^{-1},z)\to (b\gamma /c,\gamma q,b\gamma^2,x) \leqno (6.3)
$$
and renormalize (see also [20]).

If $A=q$ or $B=q$ then (6.1) becomes the recurrence for monic 
Al-Salam-Chihara polynomials [2]. For other references to this case
see [17, p. 63].

We record here the solutions to (6.1) based on (6.2) and the
solutions we obtained for $X_n$ in Sections 2 and 4.

Using (2.8) we have
$$
\leqalignno{
Q_n^{(1),\pm} (z; A,B,\delta ) \ & = \ \lim\limits_{C \to 0} C^n X_n^{(1),\pm}
\left({z\over C}; A,B,C,\delta C\right) & (6.4)\cr
& = \ {(A,B)_n\over (AB \Lambda_\pm )_n} (\Lambda_\pm )^n {}_2\phi_1
\left(\matrix{B\Lambda_\pm, Bq^n\cr AB\Lambda_\pm q^n\cr} ;
A\delta \Lambda_\pm \right),\cr
\Lambda_\pm \ & = \ \frac12 (z\pm \sqrt{z^2-\gamma^2} ) , & (6.5)\cr
\gamma \ & = \ 2(q/AB\delta )^\frac12 .\cr}
$$
Also from (2.14), (2.15) and (2.19) respectively we similarly
obtain
$$
\leqalignno{
Q_n^{(2)} (z;A,B,\delta ) &\ = \ (B)_n B^{-n} {}_2\phi_1 
\left( \matrix{B\Lambda_+, B\Lambda_-\cr q/\delta\cr};
q^{-n+1} /B \right), & (6.6)\cr
Q_n^{(3)} (z;A,B,\delta ) &\ = \ (B)_n (\delta B)^{-n} {}_2\phi_1 
\left( \matrix{B\delta \Lambda_+, B\delta \Lambda_-\cr q\delta\cr};
q^{-n+1} /B \right), & (6.7)\cr
Q_n^{(4)} (z;A,B,\delta ) &\ = \ (-q/AB\delta )^n q^{-n(n-1)/2} & (6.8)\cr
& \times (AB\delta \Lambda_+ /q, AB\delta \Lambda_- /q)_n \ {}_2\phi_2
\left( \matrix{q^{-n+1}/A, q^{-n+1}/B\cr \Lambda_+ q^{-n+1},
\Lambda_- q^{-n+1}\cr}; q/\delta \right), \cr}
$$
with $Q_n^{(4)} $ proportional to $Q_n^{(2)}$ via (5.28) and
(5.29).

Using (4.12) we obtain the explicit polynomial formula (first make
the interchange $B \ \leftrightarrow \ D$)
$$
\leqalignno{
Q_n (z;A,B,\delta ) \ = \ & \lim\limits_{C\to 0} C^n P_n (z/C;
A,\delta C, C, B) & (6.9)\cr
\ = \ & (\gamma u/2)^n {(A,B)_n\over (q)_n} \cr
& \quad \times \left\{ \sum_{\ell = 0} ^n {(q^{-n}, 2u/\gamma 
\delta, 2u /\gamma )_\ell\over (A,B)_\ell} (-1)^\ell u^{-2\ell }
q^{n\ell} q^{-\ell (\ell -1)/2}\right.\cr
& \quad \left.\times \sum_{j=0}^\ell {(A/q, B/q)_j (-1)^j u^{2j} q^{j(j+1)/2}
\over (q,2u/\gamma \delta, 2u/\gamma )_j } \right\},\cr
z \ = \ & \gamma (u+u^{-1})/2,\qquad \gamma \ = \ 2(q/AB\delta )^\frac12.\cr}
$$
With the choice of square root branch chosen so that $| \Lambda_-
/\Lambda_+ | < 1$ for $x = z/\gamma \in \bC \setminus [-1,1]$,
we have the minimal solution to (6.1) given by (6.5).  As a 
consequence we may give an explicit expression for the corresponding
continued fraction and the absolutely continuous component of the measure 
which gives its representation as a Stieltjes transform.
The calculations proceed as in Section 3 and yield the following.
For $z /\gamma = x \in \bC \setminus [-1,1]$ and $|\Lambda_-/\Lambda_+|
<1$,
\vfill\eject
$$
\leqalignno{
\qquad \quad {1\over z-a_0}{{}\atop -} {b_1^2\over z-a_1}{{}\atop -} {b_2^2\over z-a_2}{{}\atop - \cdots}
\ = \ &{AB\delta \Lambda_-\over q(1-AB \Lambda _- /q)}
{{}_2\phi_1 \left(\matrix{B\Lambda_-,B\cr
AB\Lambda_-\cr}; A\delta \Lambda_- \right)\over
{}_2\phi_1 \left(\matrix{B\Lambda_-,B/q\cr
AB\Lambda_-/q\cr}; A\delta \Lambda_- \right)} & (6.10) \cr
\ = \ & \int_{-1}^1 {\omega (t) dt\over z-\gamma t} + \
\hbox{possible pole terms} 
,\cr}
$$
with
$$
\leqalignno{
\omega(x) \ = \ & {1\over 2\pi \sqrt{1-x^2}} {(A,B, 1/u^2,u^2)_\infty\over
(A\delta \gamma u/2, A\delta \gamma/2u, AB\gamma u/2q, AB\gamma /2qu)_\infty} & (6.11)\cr
\times & \left[{}_2\phi_1 \left(\matrix{B\Lambda_- , B/q\cr
AB\Lambda_-/q\cr} ; A\delta \Lambda_- \right)
 2\phi_1 \left(\matrix{B\Lambda_+ , B/q\cr
AB\Lambda_+/q\cr} ; A\delta \Lambda_+ \right)\right]^{-1}.\cr}
$$
Note that (6.11) agrees with the weight function derived in
[2, (3.64) with $a = q(1+\delta^{-1} )/B, b=q^2/B^2\delta,c =q/AB\delta
,\gamma =B/q$ and $x^2$ replaced by $x^2/4c$].
\bigskip

\noindent{\bf 6.2  Associated Al-Salam-Carlitz I.}
\medskip

We take the $B \to \infty $ limit of (6.1) and its corresponding
solutions to obtain the case of associated Al-Salam-Carlitz I.
The recurrence becomes 
$$
R_{n+1} - [z - (1+\delta ^{-1})q^n ] R_n - {q^n\over A\delta}
(1-Aq^{n-1} ) R_{n-1} \ = \ 0 \leqno(6.12)
$$
with solutions from (6.5) -- (6.8) given by
$$
\leqalignno{
R_n^{(1)} (z;A,\delta ) \ = \ & \lim\limits_{B\to \infty} Q_n^{(1),-}
(z;A,B,\delta) & (6.13) \cr
\ = \ & {(A)_n\over (q/\delta z)_n} (-q/A\delta z)^n q^{n(n-1)/2}
{}_1 \phi_1 \left(\matrix{q/Az\delta\cr q^{n+1}/z\delta\cr};
q^{n+1}/z \right),\cr
R_n^{(2)} (z;A,\delta ) \ = \ & \lim\limits_{B\to \infty} Q_n^{(2)} & (6.14)\cr
\ = \ & (-1)^n  q^{n(n-1)/2} 
{}_1 \phi_1 \left(\matrix{q/Az\delta\cr q/\delta\cr};
zq^{-n+1} \right),\cr
R_n^{(3)} (z;A,\delta ) \ = \ & \lim\limits_{B\to \infty} Q_n^{(3)} & (6.15)\cr
\ = \ & (-\delta)^{-n}  q^{n(n-1)/2}
{}_1 \phi_1 \left(\matrix{q/Az\cr q \delta\cr};
\delta zq^{-n+1} \right),\cr
R_n^{(4)} (z;A,\delta ) \ = \ & \lim\limits_{B\to \infty} Q_n^{(4)} & (6.16)\cr
\ = \ &  (1/z)_n\ z^n
{}_1 \phi_1 \left(\matrix{q^{-n+1}/A\cr zq^{-n+1}\cr};
q/\delta \right),\cr}
$$
with $R_n^{(4)}$ proportional to $R_n^{(2)}$ via (5.29) and (5.30);
namely
$$
{}_1\phi_1 \left( \matrix{c/b\cr c\cr} ; bz \right) \ = \
{(bz)_\infty\over (c)_\infty} {}_1\phi_1 \left(\matrix{z\cr
bz\cr}; c\right). \leqno(6.17)
$$
The minimal solution to (6.12) is given by (6.13).
Using Pincherle's theorem we obtain the continued fraction
representation
$$
\leqalignno{
{1\over z-(1+\delta^{-1})} & {{}\atop +}  {q(1-A)/A\delta\over z-
(1+\delta^{-1})q} {{}\atop + } {q^2 (1-Aq)/A\delta\over
z-(1+ \delta^{-1})q^2} {{}\atop + \cdots } & (6.18)\cr
& = \ {1\over z(1-1/\delta z)} {{}_1\phi_1 \left( \matrix{q/Az\delta\cr
q/z\delta\cr} ; q/z\right)\over {}_1\phi_1 \left( \matrix{q/Az\delta\cr
1/z\delta\cr}; 1/z\right)}.\cr}
$$
This is a positive definite $ J$-fraction in the case when $0<q<1$ and
$A<1, A\delta <0$. We may then deduce that the zeros of the
${}_1\phi_1$'s in the numerator and denominator on the right side of
(6.18) are real and simple and interlace. See for example the Corollary in Section 5.4.

When $A = q$, the pole singularities in (6.18) become
explicit, since we may then use (6.17) to obtain
$$
\leqalignno{
\qquad {1\over z-(1+\delta^{-1})} & {{}\atop +}  {(1-q)/\delta\over z-
(1+\delta^{-1})q} {{}\atop + } {q (1-q^2)/\delta\over
z-(1+ \delta^{-1})q^2} {{}\atop + \cdots } & (6.18')\cr
& = \ {(q/\delta z)_\infty \over z(1/z)_\infty (1/z\delta )_\infty}
{}_1\phi_1 \left( \matrix{1/z\delta\cr
q/z\delta\cr} ; q/z\right)\cr
& = \ \sum_{n=0}^\infty \left\{ {q^n\over (z-q^n)(q)_n (\delta q)_n
(\delta^{-1} )_\infty} + {q^n\over (z-q^n/\delta )(q)_n (q/\delta )_n
(\delta )_\infty} \right\}.\cr}
$$
In the last equality we have assumed that $\delta\ne q^{-m}, m$ an
integer.
The explicit polynomial solution to (6.12) can be obtained from the
$B \to \infty$ limit of (6.9).  We get
$$
\leqalignno{
R_n (z;A,\delta ) & = \ \lim\limits_{B\to \infty} Q_n (z;A,B,\delta)
& (6.19)\cr
& = \ (-q/A\delta z)^n {(A)_n\over (q)_n} q^{n(n-1)/2} \times
\left\{ \sum_{\ell = 0}^n {(q^{-n}, 1/z\delta , 1/z)_\ell\over
(A)_\ell} \right.\cr
& \qquad \left.\times q^{-\ell (\ell -1)} (A\delta z^2 /q)^\ell
q^{n\ell } \sum_{j=0}^\ell {(A/q)_j q^{j^2} (A\delta z^2)^{-j}\over
(q, 1/z\delta , 1/z )_j} \right\}.\cr}
$$
When $A = q$, the expression for $R_n (z;q,\delta)$ must be equal to
$R_n^{(4)} (z;q,\delta )$ given by (6.16), since they are both
monic polynomial solutions.  Equating these two expressions
we obtain
$$
\leqalignno{
R_n (z; q,\delta ) \ = \ & (-\delta z)^{-n} q^{n(n-1)/2} {}_3 \phi_0
\left(\matrix{q^{-n}, 1/z \delta, 1/z\cr -\cr}; \delta z^2 q^n\right) & (6.20)\cr
\ = \ & z^n (1/z)_n {}_1\phi_1 \left( \matrix{q^{-n}\cr zq^{-n+1}\cr};
q/\delta \right).\cr}
$$
The above connection between a terminating ${}_3\phi_0$ and a terminating
${}_1\phi_1$ appears to be new.  Note that (6.20) differs also
from the standard expression in [17, (3.24.1) with $a = \delta^{-1}$
and $x$ replaced by $z$] which gives
$$
\leqalignno{
R_n (z;q,\delta ) \ = \ & (-\delta)^{-n} q^{n(n-1)/2} {}_2\phi_1
\left(\matrix{q^{-n}, z^{-1}\cr 0\cr}; qz\delta \right) & (6.21)\cr
\ = \ & z^n (1/z\delta )_n \ {}_1\phi_1 \left(\matrix{
q^{-n}\cr q^{-n+1} z\delta\cr}; q\delta \right).\cr}
$$
For the last equality we have used (5.29) with $c = 0$.
\medskip

\noindent{\bf 6.3 Limit Al-Salam-Carlitz I.}
\medskip

We now take the $A \to \infty$ limit of (6.12) and its corresponding
solutions.  The recurrence becomes
$$
S_{n+1} - [z - (1+\delta^{-1}) q^n ] S_n + {q^{2n-1}\over \delta}
S_{n-1} \ = \ 0 \leqno(6.22)
$$
with solutions
$$
\leqalignno{
S_n^{(1)} (z; \delta ) \ = \ & \lim\limits_{A \to \infty} R_n^{(1)}
(z;A,\delta ) & (6.23)\cr
\ = \ & {q^{n^2}\over (q/\delta z)_n (\delta z)^n} {}_1\phi_1
\left(\matrix{0\cr q^{n+1}/z\delta\cr}; q^{n+1}/z\right),\cr
S_n^{(2)} (z; \delta ) \ = \ & \lim\limits_{A \to \infty} R_n^{(2)}
(z;A,\delta ) & (6.24)\cr
\ = \ & (-1)^n q^{n(n-1)/2} {}_1\phi_1
\left(\matrix{0\cr q/\delta\cr}; zq^{-n+1}\right),\cr
S_n^{(3)} (z; \delta ) \ = \ & \lim\limits_{A \to \infty} R_n^{(3)}
(z;A,\delta ) & (6.25)\cr
\ = \ & (-\delta)^{-n} q^{n(n-1)/2} {}_1\phi_1
\left(\matrix{0\cr q\delta\cr}; \delta zq^{-n+1}\right),\cr
S_n^{(4)} (z; \delta ) \ = \ & \lim\limits_{A \to \infty} R_n^{(4)}
(z;A,\delta ) & (6.26)\cr
\ = \ & (1/z)_n z^n {}_1\phi_1
\left(\matrix{0\cr zq^{-n+1}\cr}; q/\delta \right).\cr}
$$
Note that $S_n^{(4)}$ is proportional to $S_n^{(2)}$ via the
transformation (6.17) which yields the identity
$$
{}_1\phi_1 \left( {0\atop c}; z\right) \ = \ {(z)_\infty\over
(c)_\infty} {}_1\phi_1 \left( {0\atop z} ; c\right).\leqno(6.27)
$$

The minimal solution to (6.22) is $S_n^{(1)} (z;\delta )$.
Using Pincherle's theorem we then obtain the evaluation of the 
continued fraction associated with (6.22). Namely

$$
\leqalignno{
{1\over z-(1+\delta^{-1})}& {{}\atop -}\ {q/\delta\over z-(1+\delta^{-1}
)q }\ {{}\atop -}\ {q^3/\delta\over z-(1+\delta^{-1} )q^2} \ {{}\atop
- \cdots } & (6.28)\cr
& = \ {1\over z(1-1/\delta z)} {{}_1\phi_1 \left(\matrix{0\cr
q/z\delta\cr};q/z\right)\over {}_1\phi_1 \left(\matrix{0\cr
1/z\delta\cr} ; q/z\right)}.\cr}
$$
This is a positive definite $J$-fraction if $0<q<1$ and $\delta >0$ 
and we may then deduce that the zeros of the ${}_1\phi_1$'s
on the right of (6.28) are real, simple and interlacing. See the Corollary in Section 5.4.

The monic polynomials are given explicitly by
$$
\leqalignno{
S_n (z;\delta ) \ = \ & \lim\limits_{A\to \infty} R_n (z;A,\delta) & (6.29)\cr
\ = \ & {(z\delta)^{-n} q^{n^2}\over (q)_n} \left\{ \sum_{\ell = 0}^n
(q^{-n}, 1/z\delta , 1/z)_\ell (-\delta)^\ell q^{-3\ell (\ell -1)/2}
\right.\cr
& \quad \left. \times z^{2\ell} q^{\ell (n-1)} \sum_{j=0}^\ell 
{z^{-2j} (-\delta)^{-j} q^{3j(j-1)/2}\over (q, 1/z\delta,
1/z)_j} \right\} .\cr}
$$
A simpler expression is obtained if one applies the generating
function method of Section 4 directly to (6.22).  This results in
$$
\leqalignno{
S_n(z;\delta)\ = \ & {\delta^{-n} q^{n(n+1)/2}\over (q)_n} \left\{ \sum_{\ell = 0}
^n (q^{-n})_\ell (1/z)_\ell (-\delta z)^\ell q^{-\ell (\ell -1)/2}\right. & (6.30) \cr
& \quad \left.\times \sum_{j=0}^\ell {(-z\delta )^{-j} q^{j(j-1)/2}\over
(1/z)_j (q)_j} \right\} .\cr}
$$
\medskip

\noindent{\bf 6.4 Associated continuous $q$-Hermite.}
\medskip

If we multiply (6.1) by $B^\frac12$, replace $zB^\frac12$ by $z$,
renormalize and let $B\to 0$, we get the recurrence
$$
H_{n+1} - zH_n + {q\over A\delta} (1-Aq^{n-1} ) H_{n-1} \ = \ 0.
\leqno(6.31)
$$
This will become the continuous $q$-Hermite case if we
set $A = q$.  The solutions to (6.31) are obtained as limits of the
solutions to (6.1).  They are
$$
\leqalignno{
H_n^{(1),\pm} (z;A,\delta ) \ = \ & \lim\limits_{B\to 0} B^{n/2}
Q_n^{(1),\pm} (z/B^\frac12; A,B,\delta ) & (6.32)\cr
\ = \ & (A)_n (\mu_\pm )^n {}_1\phi_1 \left(\matrix{Aq^n\cr 0\cr};
A\delta \mu_\pm^2\right),\cr
H_n^{(2)} (z;A,\delta ) \ = \ & (\mu_- )^n {}_2\phi_0
\left(\matrix{q^{-n+1}/A,0\cr - \cr}; q^n/\delta \mu_-^2\right),
& (6.33)\cr}
$$
where
$$
\mu_\pm \ = \ (z\pm \sqrt{z^2-\gamma^{'2})}/2,\quad \gamma'\ = \
2\sqrt{q/A\delta}.\leqno(6.34)
$$
To derive (6.32) we first made the transformation (5.28) before taking
the limit and the transformation (5.30) in the form
$$
{}_0\phi_1 \left({-\atop c}; cz\right) \ = \ {1\over (c)_\infty}
{}_1\phi_1 \left({z\atop 0}; c\right) \leqno(6.35)
$$
after the limit was taken.

To derive (6.33) we transformed either $Q^{(2)}_n$ or $Q^{(3)}_n$
using an iterate of (5.28), discarded factors which were $n$-independent,
multiplied by $B^{n/2}$ and then let $B \to 0$.  Note that (6.33)
is only a formal solution unless it terminates by having $A = q$.

For the general polynomial solution we take the limit of (6.9) to obtain
$$
\leqalignno{
H_n(z;A,\delta ) \ = \ & \lim\limits_{B\to 0} B^{n/2} Q_n (z;A,B,
\delta) & (6.36) \cr
\ = \ & (\gamma' u/2)^n {(A)_n\over (q)_n} \left\{ \sum_{\ell = 0}^n
{(q^{-n})_\ell\over (A)_\ell} (-1)^\ell u^{-2\ell} q^{n\ell}\right.\cr
& \quad \left.\times q^{-\ell (\ell -1)/2} \sum_{j=0}^\ell {(A/q)_j (-1)^j u^{2j}
q^{j(j+1)/2} \over (q)_j}\right\} , \cr
z \ = \ & {\gamma'\over 2} (u+u^{-1}),\quad \gamma' \ = \ 2(q/A\delta )^\frac12.\cr}
$$
Note that when $A = q$, the expressions (6.33) and (6.36) become
equal.

For $z = \gamma' x,\ x \in \bC \setminus [-1,1]$ and with the 
square root branch chosen so that $| \mu_- / \mu_+ | < 1$,
we have $H_n^{(1),-} (z; A,\delta )$ as a minimal solution to
(6.31) and hence the continued fraction representation
$$
\leqalignno{
& {1\over z} {{}\atop -} {q(1-A)/A\delta\over z} {{}\atop -} 
{q(1-Aq)/A\delta\over z} {{}\atop -} {q(1-Aq^2)/A\delta\over z}
{{}\atop - \cdots } & (6.37)\cr
& = \ {A\delta \mu_-\over q} {{}_1\phi_1 \left(\matrix{A \cr 0 \cr}; A\delta
\mu_-^2 \right)\over {}_1\phi_1 \left( \matrix{A/q\cr 0\cr};
A\delta \mu_-^2 \right)}\cr
& = \ \int_{-1}^1 {\omega (t) dt\over z - \gamma' t} +
\hbox{possible pole terms.}\cr}
$$

Repeating the method of Section 3, we find that
$$
\leqalignno{
\omega (x) \ = \ & {1\over 2\pi \sqrt{1-x^2}} {(A,u^2,u^{-2})_\infty
\over {}_1\phi_1 \left(\matrix{A/q\cr 0\cr}; u^2/q\right) {}_1\phi_1
\left(\matrix{A/q\cr 0\cr}; 1/u^2q\right)}, & (6.38)\cr
z/\gamma' \ = \ & x \ = \ \frac12 (u+u^{-1} ).\cr}
$$
When $A = q$ the denominator ${}_1\phi_1$'s in (6.37)
and (6.38) become equal to 1 and there are no pole terms. This is then the
continuous $q$-Hermite case.  See [17, p.88] for a list of references.
\bigskip

\noindent{\bf 6.5 Limit $q$-Hermite.}
\medskip

We take the $A \to \infty$ limit of the results in the previous Section
6.4.  This gives the recurrence relation
$$
T_{n+1} - zT_n - {q^n\over \delta} T_{n-1} \ = \ 0 \leqno(6.39)
$$
with minimal solution
$$
\leqalignno{
T_n^{(1)} (z;\delta ) \ = \ & \lim\limits_{A\to \infty} H_n^{(1),-} (z;A,\delta ) & (6.40) \cr
\ = \ & (-1)^n q^{n(n-1)/2} (q/\delta z)^n {}_0\phi_1 \left(
{-\atop 0}; {q^{n+2}/\delta z^2}\right)\cr}
$$
and the polynomial solution
$$
\leqalignno{
T_n (z;\delta ) \ = \ & \lim\limits_{A\to \infty} H_n (z;A,\delta ) & (6.41) \cr
\ = \ & {(-z)^{-n} q^{n(n-1)/2} \over (q)_n} (q/\delta)^n \left\{
\sum_{\ell = 0}^n (q^{-n})_\ell \right.\cr
& \left.\qquad \times z^{2\ell} q^{n\ell} (\delta/q)^{\ell} q^{-\ell (\ell -1)}
\sum_{j=0}^\ell {z^{-2j} \delta^{-j} q^{j^2}\over (q)_j}\right\}.\cr}
$$
Using the minimal solution we get the continued fraction result
$$
{1\over z} {{}\atop +} {q/\delta\over z} {{}\atop +} {q^2/\delta
\over z} {{}\atop + \cdots } 
 \ = \  {1\over z} {{}_0\phi_1 \left( {-\atop 0}; q^2 /\delta z^2 \right)\over
{}_0\phi_1 \left( {-\atop 0}; q/\delta z^2\right)}. \leqno(6.32)
$$
This is a positive definite $J$-fraction if $0<q<1$ and $\delta <0$ 
and one may then has simple real interlacing zeros for the ${}_0\phi_1$'s
on the right of (6.32) (see the Corollary in Section 5.4). 
\bigskip

\noindent{\bf 6.6 Associated continuous big $q$-Hermite.}
\bigskip

In (6.1) let $\delta = 1 /aB$ and then let $B \to 0$ to get
the recurrence
$$
C_{n+1} - (z-q^n) C_n + {aq\over A} (1-Aq^{n-1}) C_{n-1} \ = \ 0.
\leqno(6.43)
$$
The solutions to (6.43) may be obtained from the solutions of
(6.1) by using
$$
C_n(z;A,a) \ = \ \lim_
{B\to 0} Q_n (z;A,B,1/aB). \leqno(6.44)
$$
First making the parameter interchange $A \leftrightarrow B$
in (6.5), then writing $\delta = 1/aB$ and taking limit 
$B \to 0$ we obtain the solution
$$
\leqalignno{
C_n^{(1),\pm} (z;A,a) \ = \ & \lim\limits_{B\to 0} {(A,B)_n\over
(AB\Lambda_\pm' )_n} (\Lambda'_\pm )^n {}_2\phi_1 \left(\matrix{
A\Lambda_\pm' , Aq^n\cr AB\Lambda_\pm' q^n\cr}; {\Lambda_{\pm}'\over
a}\right), & (6.45)\cr
\ = \ & (A)_n (\Lambda'_\pm )^n {}_2\phi_1 \left(\matrix{A\Lambda_\pm',
Aq^n\cr 0\cr}; {\Lambda_\pm'\over a}\right),\cr
\Lambda_\pm' \ = \ & \frac12 (z\pm \sqrt{z^2-\gamma_1^2})\cr
\gamma_1 \ = \ & 2\sqrt{{aq\over A}}.\cr}
$$
Also from (6.2) we similarly obtain
$$
\leqalignno{
C_n^{(2)} (z;A,a) \ = \ & \lim\limits_{B\to 0} Q_n^{(2)} (z; B,A,
1 /aB) & (6.46)\cr
\ = \ & {\rm const}\ (\Lambda_+' )^n (A\Lambda'_- /aq)_n\ 
{}_2\phi_1 \left(\matrix{q^{-n+1}/A,0\cr \Lambda_+' q^{-n+1}\cr}
; A\Lambda_-' \right)\cr}
$$
which is a polynomial solution when $A = q$.

In (6.7) we first apply the transformation [5,(III.2), p. 241], write
$\delta = 1/aB,\ B = q^m$ and let $m\to \infty$.  We obtain the solution
$$
C_n^{(3)} (z;A,a) \ = \ {\rm const}\ (\Lambda_+' )^n {}_2\phi_0
\left(\matrix{q^{-n+1}/A,\Lambda_+'/a\cr -\cr};
{A^2\Lambda_-^{'2} q^{n-2}\over a}\right), \leqno(6.47)
$$
which is again a polynomial solution for $A = q$.
Using (6.8), transformation (5.30) and taking limits we have the
solution
$$
C_n^{(4)} (z; A,a) \ = \ {\rm const}\ (\Lambda_+' )^n (A \Lambda'_- /aq)_n
\times {}_2\phi_1 \left(\matrix{q^{-n+1}/A,0\cr \Lambda_+' q^{-n+1}\cr};
A \Lambda_-' \right). \leqno(6.48)
$$
Also from (6.9) we obtain the explicit polynomial solution
$$
\leqalignno{
C_n(z;A,a) \ = \ & \lim\limits_{B\to 0} Q_n (z;A,B,{1\over aB} ) & (6.49)\cr
\ = \ & (\gamma_1 u/2)^n {(A)_n\over (q)_n} \left\{ \sum_{\ell=0}^n
{(q^{_n}, 2u/\gamma_1)_\ell\over (A)_\ell} (-1)^\ell u^{-2\ell}
q^{n\ell} q^{-\ell (\ell -1)/2}\right.\cr
& \quad \times \left.\sum_{j=0}^\ell {(A/q)_j (-1)^j u^{2j} q^{j(j+1)/2}\over
(q,2u / \gamma_1)_j}\right\},\cr
\gamma_1 \ = \ & 2\sqrt{aq /A}.\cr}
$$

If in (6.43) we write $z = x\gamma_1 = 2(aq /A)^\frac12 x$, and renormalize,
we obtain the recurrence for the associated continuous big 
$q$-Hermite polynomials (see [17, (3.18.4)]) viz., the
relation
$$
\xi_{n+1} - (2x - bq^n )\xi_n + (1-Aq^{n-1}) \xi_{n-1} \ = \ 0 
\leqno(6.50)
$$
with $b = (A/aq)^\frac12$.  Writing $A = q$ in the renormalized
solution given by (6.47) yields continuous big $q$-Hermite
polynomials [17, (3.18.1)].

Choosing, for $x = z/\gamma_1 \in \bC \setminus [-1,1]$, the
square root branch for which $|\Lambda_-' /\Lambda_+' | < 1$,
the minimal solution of (6.43) is given by (6.45).
Consequently we have the continued fraction representation
$$
\leqalignno{
{1\over z-a_0} {{}\atop -} {b_1\over z-a_1} {{}\atop -} {b_2\over
z-a_2} {{}\atop -} {{}\atop \cdots} \ = \ &
{A\Lambda_-'\over aq} {{}_2\phi_1 \left(\matrix{A,A\Lambda_-'\cr
0\cr}; \Lambda_-' /a \right)\over {}_2\phi_1 \left(\matrix{A,A\Lambda_-'\cr
0\cr}; \Lambda_-' /a\right)} & (6.51)\cr
\ = \ & \int_{-1}^1 {\omega (t) dt\over z-\gamma_1t} +
\hbox{possible pole terms},\cr}
$$
with
$$
\omega(x)  =  {1\over 2\pi \sqrt{1-x^2}} {(A)_\infty (1/u^2,u^2)_\infty\over
(\gamma_1 u/2a, \gamma_1 /2au)_\infty }\left[{}_2\phi_1 \left(\matrix{A/q,
A\Lambda_-'\cr 0\cr} ;{\Lambda_-'\over a}\right)
{}_2\phi_1 \left(\matrix{ A/q,A\Lambda_+'\cr 0\cr}
;{\Lambda_+'\over a} \right)\right]^{-1} \leqno(6.52)
$$
where
$$
\eqalign{
x \ = \ & \frac12 (u+u^{-1} )\cr
\gamma_1 \ = \ & 2(aq/A )^\frac12\cr
a_n \ = \ & q^n\cr
b_n^2 \ = \ & {qa\over A} (1-Aq^{n-1}).\cr}
$$

In the particular case $A = q$,
$$
\omega(x) \ = \ {1\over 2\pi \sqrt{1-x^2}} {(q)_\infty (1/u^2,u^2)_\infty
\over (u/\sqrt{a}, 1/u\sqrt{a})_\infty}.\leqno(6.53)
$$
This weight function agrees with the form given in [17, (3.18.2]
with $b$ in (6.50) replaced by $1/\sqrt{a}$ and taking
into consideration the normalization factor
$(A/aq)^{n/2} \ = \ a^{-n/2} = b^n$ used in (6.50) and also that
$x$ is replaced by $bx$.
\vfill\eject
\noindent{\bf 6.7  $q$-Bessel order.}
\bigskip

The limit $A \to \infty$ of the recurrence (6.43) gives
$$
B_{n+1} - (z - q^n ) B_n - aq^n B_{n-1} \ = \ 0. \leqno(6.54)
$$
It is clear that (6.54) will give real orthogonal polynomials only
for $0 < q < 1,\ a < 0$.  The solutions of (6.54) can be
obtained by taking $A \to \infty$ limits of (6.45)- (6.48).
We have
$$
\leqalignno{
B_n^{(1)} (z;a) \ = \ & \lim\limits_{A\to \infty} C_n^{(1),-}
(z;A,a) \ = \ (-aq/z)^n q^{n(n-1)/2} {}_1\phi_1
\left(\matrix{aq/z\cr 0\cr}; q^{n+1}/z\right) & (6.55)\cr
\ = \ & (-aq/z)^n q^{n(n-1)/2} (q^{n+1}/z)_\infty \ 
{}_0\phi_1\left(\matrix{- \cr q^{n+1}/z\cr} ; aq^{n+2}/z^2\right), & (6.56)\cr}
$$
using the transformation (6.35). Also
$$
B_n^{(2)} (z;a) \ = \ \lim\limits_{A\to \infty} C_n^{(2)} 
(z;A,a) \ = \ {\rm const }\ z^n (1/z)_n {}_2\phi_1 
\left(\matrix{0,0\cr zq^{-n+1}\cr}; aq/z\right),\leqno(6.57)
$$
$$
B_n^{(3)} (z;a) \ = \ \lim\limits_{A\to\infty} C_n^{(3)} (z;A,a)
\ = \ {\rm const} \ z^n {}_2\phi_0 \left( \matrix{0,z/a\cr -\cr};
aq^n/z^2 \right),\leqno(6.58)
$$
$$
 \lim\limits_{A\to \infty} C_n^{(4)} (z; A,a)
\ = \ {\rm const}\ B_n^{(2)} (z;a).\leqno(6.59)
$$
Using (6.49) we obtain the explicit polynomial 
$$
\leqalignno{
B_n(z;a) \ = \ & \lim\limits_{A \to \infty} C_n(z;A,a) & (6.60)\cr
\ = \ & {(-a/z)^n\over (q)_n} q^{n(n+1)/2} \left\{ \sum_{\ell = 0}^n
(q^{-n},1/z)_\ell \right. \cr
& \quad \left. \times\ q^{-\ell^2} q^{n\ell} \left({z^2\over a}\right)^\ell
\sum_{j=0}^\ell {q^{j^2}(a/z^2)^j\over (q)_j (1/z)_j} \right\}.\cr}
$$
We now demonstrate the relationship of the above solutions
with Jackson's $q$-analogues of Bessel functions (see [5, Exercise 1.24,
p. 25]).  Using the notation of [5] and writing
$$
q^\nu \ = \ z^{-1} , \ - x^2/4 \ = \ aq/z,
$$
we find from (6.57) and (6.56) that
$$
J_{-\nu-n}^{(1)} (2i (aq/z)^\frac12; q) \ = \ {\rm const}\
(-1)^n q^{-n(n+4)/2} \sqrt{z} (-i \sqrt{z/a})^{n+\nu} 
{(q^2/z)_n\over (1/z)_n} B_n^{(2)} (z;a) \leqno(6.61)
$$
and
$$
J_{\nu+n}^{(2)} (2i (aq/z)^\frac12; q) \ = \ {\rm const}\
(-1)^n {q^{-n(n-1)/2}\over \sqrt{z}} (z/a\sqrt{q})^n
(i \sqrt{a/z})^{n+\nu}  B_n^{(1)} (z;a). \leqno(6.62)
$$
This shows that (6.54) is connected with the recurrence for
$q$-Bessel functions with $z$ appearing both in the
arguement and the order of the $q$-Bessel function.
Thus we choose to call this the $q$-Bessel order case.  The
continued fraction representation obtained with the help
of minimal solution of (6.54) given by (6.55) is
$$
{1\over z-1} {{}\atop +} {aq\over z-q} {{}\atop +} {aq^2\over
z-q^2} {{}\atop +} {{}\atop \cdots} \ = \
{1\over (z-1)} {{}_0\phi_1 \left(\matrix{-\cr q/z\cr}; aq^2/z^2 \right)\over
{}_0\phi_1 \left(\matrix{-\cr 1/z\cr}; aq/z^2 \right)}. \leqno(6.63)
$$
This is a positive definite $J$-fraction if $0<q<1$ and $a<0$ 
and the zeros of the ${}0\phi_1$'s on the right are then real simple and
interlacing.

A $q\to 1$ limiting case is due to Maki [19] (see also [21]).
\bigskip
\noindent{\bf 7. Limit Askey-Wilson.}
\bigskip

In this section we give the connection between solutions to the
associated Askey-Wilson and the associated continuous dual $q$-Hahn polynomial 
recurrence relations.  The associated Askey-Wilson polynomial
recurrence relation in monic form is given by
$$
\leqalignno{
p_{n+1} & (x) - (x - a_n) p_n (x)+ b_n^2 p_{n-1} (x) \ = \ 0, & (7.1)\cr
a_n \ = \ & -A_n - B_n + {a\over 2} + {1\over 2a},\cr
b_n^2 \ = \ & A_{n-1} B_n,\cr}
$$
$$
\eqalign{
A_n \ = \ & {(1-abcd q^{n+\alpha -1} ) (1-abq^{n+\alpha} )(1-acq^{n+\alpha})
(1-adq^{n+\alpha} )\over 2a(1-abcdq^{2n+2\alpha-1})(1-abcdq^{2n+2\alpha})},\cr
B_n \ = \ & {a(1- q^{n+\alpha } ) (1-bcq^{n+\alpha-1} )(1-bdq^{n+\alpha-1})
(1-cdq^{n+\alpha-1} )\over 2(1-abcdq^{2n+2\alpha-2})(1-abcdq^{2n+2\alpha-1})}.\cr}
$$
When $\alpha = 0$, (7.1) reduces to the non-associated monic
Askey-Wilson case [3].  If one further puts $d = 0$, one has
the monic continuous dual $q$-Hahn case [3, pp.3, 28].  Our
associated dual $q$-Hahn case may be obtained from (7.1) by
multiplying (7.1) by $k^{-1} = 2\sqrt{q/ABCD},$ replacing
$(a,b,c,d,q^{\alpha})$ by $(2kq/AD, 2kq/AC, 2kq/AB, 0,A/q)$ and
renormalizing to monic form.  Note that `$k$' is now used to
denote the `$\alpha$' of previous Sections so as to avoid
confusion with the `$\alpha$' in (7.1).

There are two papers [15], [9], which deal with the
associated Askey-Wilson polynomial case.  We now indicate
how our associated dual $q$-Hahn recurrence solutions in
Section 2 are connected with the solutions in [15] and [9].

Using the solutions to (7.1) which Ismail and Rahman have
given in [15], we have the following limiting cases.
$$
P_n (z; A,B,C,D) \ = \ \lim\limits_{d\to 0} {(B,C)_n\over
k^n} p_n(x; {2kq\over AD}, {2kq\over AC}, {2kq\over
AB},d,A/q)\leqno(7.2)
$$
where
$$
p_n(x)=p_n (x; a,b,c,d,q^{\alpha}) \ = \ {1\over (2a)^n}
p_n^{(\alpha )} (x)
$$
and $p_n^{(\alpha )} (x)$ is given explicitly in [15,
(4.15)].  The calculation yields
$$
\leqalignno{
P_n(z; A,B,C,D) \ = \ & {(B,C)_n\over (BC)^n} \left\{
\sum_{k=0}^n {(q^{-n}, \sqrt{BCq /AD} u,
\sqrt{BCq /AD}  u^{-1})_k \ q^k\over (q,B,C)_k}\right.
& (7.3)\cr
& \left. \times \sum_{j=0}^{n-k} {(A /q, D /q, q^{k+1},
q^{k-n})_j\over (q,Cq^k , Bq^k, q^{-n})_j} \left({BCq\over
AD}\right)^j \right\}.\cr}
$$
Ismail and Rahman have also obtained non-polynomial
solutions $r_{n+\alpha} (u)$ and $s_{n+\alpha} (u)$
which correspond to $p_n^{(\alpha )} (x)$.  Using [15,
(1.12), (1.13)] we obtain
$$
X_n^{(1),\pm} (z; A,B,C,D) \ = \ \lim\limits_{d\to 0}
C_{1}^{\pm}{(B,C)_n\over (BC)^n} \left( (BC\lambda_\pm)
^{-\alpha} s_{n+\alpha}(u^{\pm 1})\right)\leqno(7.4)
$$
where
$$
C_{1}^{\pm} \ = \ {(A,D,BC\lambda_\pm )_\infty\over
(BCD\lambda_\pm , ABC\lambda_\pm , AD\lambda_\pm)_{\infty}},
$$
and
$$
X_n^{(5)} (z;A,B,C,D)\ = \ \lim\limits_{d\to 0} C_5
{(B,C)_n\over (BC)^n} \left( (BC\lambda_-)^{-\alpha}
r_{n+\alpha} (u)\right) \leqno(7.5)
$$
where
$$
C_5 \ = \ {(C,AD\lambda_+, Dq \lambda_-, ABC\lambda _+
/q)_\infty\over (B,C, q/B, Dq/C )_\infty}.
$$

In [9], the parameters $(a,b,c,d,q^\alpha )$ of [15]
are replaced by $(\alpha,\beta,\gamma,\delta,\epsilon)$. The
solutions $X_n^{(r)} $ to (7.1) in [9] will now be
denoted by $\tX_n^{(r)}$.  The $\delta \to 0$ limits
of some of these solutions are given below in terms of the
solutions $X_n^{(1)}, X_n^{(2)}, X_n^{(4)}, X_n^{(5)}$ and
$X_n^{(6)}$ of Section 2.  We have
$$
\lim\limits_{\delta \to 0} \tX_n^{(6)} (u^{\pm 1}) \ = \
D_{6}^{\pm}k^n X_n^{(1),\pm} (z;A,B,C,D),\leqno(7.6)
$$
$$
D_{6}^{\pm} \ = \ {(BCD\lambda_\pm , AD\lambda_\pm,
ABC\lambda_\pm)_\infty\over (A,B,C,D,q/u^{\pm 2})_\infty}
$$
$$
\leqalignno{
\lim\limits_{\delta\to 0} \tX_n^{(5)} (u) \ = \ & D_5
k^nX_n^{(5)} (z;A,B,C,D) + D'_5 k^nX_n^{(2)}
(z;A,B,C,D), & (7.7)\cr
D_5 \ = \ & {({A\over q} u^2, {q^2\over Au^2}, {q\over D}, AC\lambda_-,
AB\lambda_-,{Dq\over C}, {q\over B})_\infty\over
({ABC\over q}\lambda_+ , {ABD\over q} \lambda_+, {ACD\over q}
\lambda_+, Cq\lambda_-, Bq\lambda_-, {q\over u^2} )_\infty}\cr
& \times {1\over ({A\over D}, Dq\lambda_-, AD\lambda_+)_\infty},\cr
D'_5 \ = \ & {({A\over q} u^2, {q^2\over Au^2}, {q\over A}, Aq\lambda_-,
BD\lambda_-,CD\lambda_-, {Dq\over 2k})_\infty\over
({ABC\over q}\lambda_+ , {ABD\over q} \lambda_+, {ACD\over q}
\lambda_+, Cq\lambda_-, Bq\lambda_-, {q\over u^2} )_\infty}\cr
& \times {({q^2\over Du}, {q\over B}, {Aq\over C},{2k\over D},{Du\over q} )_\infty\over
( Dq\lambda_-, {D\over A},
{q^2\over Au}, {Aq\over 2k}, {Au\over q}, {2k\over A}, Aq\lambda_-,
AD\lambda_+)_\infty}.\cr}
$$
$$
\leqalignno{
\lim\limits_{\delta\to 0} \tX_n^{(2)} (u) \ = \ & D_2 k^nX_n^{(6)} (z; C,B,A,D), & (7.8) \cr
D_2 \ = \ &{2\over u} {(A^2 q\lambda_+ , Aq\lambda_- )_\infty\over
({Aq\over B}, {Aq\over C}, {Aq\over D}, {BCD\lambda_-\over q} )_\infty}.\cr}
$$
$$
\leqalignno{
\lim\limits_{\delta\to 0} \tX_n^{(1)} (u) \ = \ & D_1
k^nX_n^{(5)} (z;A,B,C,D) + D'_1 k^nX_n^{(1)}
(z;A,B,C,D), & (7.9)\cr
D_1 \ = \ & {(q /B ,BCDq\lambda_+/ A, BC\lambda_-,
Dq /C)_\infty\over
(D,Cq /A, Bq /A, Dq \lambda_+, 
ABC \lambda_- /q, AD\lambda_- )_\infty}\cr
D'_1 \ = \ & {({q\over A} ,{BCq\over 2k}, {2k\over BC},
{Dq\over 2k}, BD\lambda_+, CD\lambda_+, {BCDq\over A}\lambda_+,
BCD\lambda_-,{2k\over D})_\infty\over
(B,C,D,{Dq\over A}, {Cq\over A}, {Bq\over A}, u, {q\over u},
{uq^2\over A} , {A\over uq},D\lambda_+ )_\infty}.\cr}
$$

Parameter interchanges in [9] in some of the cases yield new
solutions of (7.1).  For example an interchange of $\delta
\leftrightarrow \beta $ in solution $\tX_n^{(1)} (u)$ gives a new
solution $\appX_n^{(1)} (u)$ for which we have
$$
\leqalignno{
\lim\limits_{\delta\to 0} \appX_n^{(1)} (u) \ = \ & D''_1
k^nX_n^{(4)} (z;A,B,C,D) , & (7.10) \cr
D_1'' \ = \ & {(q /B, Cq /D )_\infty\over
(C,Cq\lambda_+, AC\lambda_-, ABD\lambda_- /q )_\infty}.\cr}
$$
If we instead make the interchange $\alpha \leftrightarrow \delta$
in solution $\tX_n^{(1)} (u)$, then the $\delta \to 0$ limit would
be proportional to $k^nX_n^{(5)} (z;A,B,C,D)$.

\vfill\eject
\centerline{\bf REFERENCES}
\medskip
\frenchspacing

\item{1.\ } R. Askey, Continuous Hahn polynomials, {\it J. Phys. A.
} {\bf 18} (1985), L1017-L1019.

\item{2.\ } R. Askey and M.E.H. Ismail, Recurrence relations, continued
fractions and orthogonal polynomials, {\it Mem. Amer. Math. Soc.}, {\bf 300} (1984).

\item{3.\ } R. Askey and J. Wilson, Some basic hypergeometric 
orthogonal polynomials, {\it Mem. Amer. Math. Soc.},  {\bf 319}
 (1985).

\item{4.\ } T.S. Chihara, An Introduction to Orthogonal Polynomials,
Gordon and Breach, New York, 1978.

\item{5.\ } G. Gasper and M. Rahman, Basic Hypergeometric Series,
Cambridge University Press, Cambridge, 1990.

\item{6.\ } W. Gautschi, Computational aspects of three-term
recurrence relations, {\it SIAM Review} {\bf 9} (1967), 24-82.

\item{7.\ } D.P. Gupta, M.E.H. Ismail and D.R. Masson,
Associated continuous Hahn Polynomials, {\it Canadian J. of Math.}
{\bf 43} (1991), 1263-1280.

\item{8.\ } D.P. Gupta, M.E.H. Ismail and D.R. Masson,
Contiguous relations, basic hypergeometric functions and
orthogonal polynomials II. Associated big $q$-Jacobi
polynomials, {\it J. Math. Anal. Appl.} {\bf 171} (1992),
477-497.

\item{9.\ } D.P. Gupta and D.R. Masson, Solutions to the associated
q-Askey-Wilson polynomial recurrence relation in: {\it Approximation
and Computation} {\bf 158}, Birkhauser, Cambridge, 1994, to appear.
\item{10.\ } W. Hahn, \"Uber Orthogonalpolynome, die $q$-Differenzengleichungen
gen\"ugen, {\it Math. Nachr.} {\bf 2} (1949), 4-34.

\item{11.\ } W. Hahn, \"Uber Polynome, die gleichzeitig zwei verschiedenen
Orthogonalsystemen angeh\"oren, {\it Math. Nachr.} {\bf 2} (1949), 263-278.

\item{12.\ } M.E.H. Ismail, J. Letessier, D.R. Masson
and G. Valent, Birth and death processes and orthogonal 
polynomials in: {\it Orthogonal Polynomials}, P. Nevai, ed., Kluwer
Academic Publishers, Dordrecht, 1990, 229-255.

\item{13.\ } M.E.H. Ismail, J. Letessier and G. Valent, Quadratic
birth and death models and associated continuous dual Hahn polynomials,
{\it SIAM J. Math. Anal.} {\bf 20} (1989), 727-737.

\item{14.\ } M.E.H. Ismail and C.A. Libis, Contiguous relations, basic
hypergeometric functions and orthogonal polynomials I, {\it J.
Math. Anal. Appl.} {\bf 141} (1989), 349-372.

\item{15.\ } M.E.H. Ismail and M. Rahman, The associated Askey-Wilson 
polynomials, {\it Trans. Amer. Math. Soc.} {\bf 328} (1991), 201-237.

\item{16.\ } W.B. Jones and W.J. Thron, Continued Fractions:
Analytic Theory and Applications, Addison-Wesley, Reading,
Mass., 1980.

\item{17.\ } R. Koekoek and R.R. Swarttouw, The Askey-scheme of 
hypergeometric orthogonal polynomials and its $q$-analogue,
{\it Reports of the Faculty of Technical Mathematics and
Informatics} no. 94-05, Delft, 1994.

\item{18.\ } C.A. Libis, Private communication.

\item{19.\ } D. Maki, On constructing distribution functions with 
applications to Lommel polynomials and Bessel functions, {\it Trans. 
Amer. Math. Soc.} {\bf 130} (1968), 281-297.

\item{20.\ } D.R. Masson, Difference equations revisited in:
{\it Canadian Mathematics Society Conference Proceedings} Vol. 9,
J.S. Feldman and L.M. Rosen, eds., AMS, Providence RI, 1988, 73-82.

\item{21.\ } \underbar{\hskip .5in}, Difference equations, continued fractions,
Jacobi matrices and orthogonal polynomials in: {\it Non-linear
Numerical Methods and Rational Approximation}, 1987, A. Cuyt,
ed., D. Reidel, Dordrecht, 1988, 239-257.

\item{22.\ } \underbar{\hskip .5in}, Associated Wilson polynomials,
{\it Constructive Approximation} {\bf 7} (1991), 521-534.

\item{23.\ } \underbar{\hskip .5in}, Wilson polynomials and some continued
fractions of Ramanujan, {\it Rocky Mountain J. Math.}
{\bf 21} (1991), 489-499.

\item{24.\ } H. S. Wall, Analytic Theory of Continued Fractions,
D. Van Nostrad, Princeton, NJ, 1948.

\item{25.\ } J. Wilson, Hypergeometric series, recurrence 
relations and some orthogonal functions, Doctoral Dissertation,
University of Wisconsin, Madison, 1978.
\end